\documentclass{article}

\usepackage{amssymb,amsmath,xcolor,graphicx,xspace,colortbl,rotating} 
\usepackage{boxedminipage}  
\usepackage{graphics}  
\usepackage{ragged2e}  
\usepackage{tabulary}  
\usepackage{varioref}  
\usepackage{xcolor}  
\graphicspath{{NilpotentBorbitsAn_graphics/}{NilpotentBorbitsAn_tcache/}{NilpotentBorbitsAn_gcache/}}
\DeclareGraphicsExtensions{.pdf,.eps,.ps,.png,.jpg,.jpeg}\newtheorem {theorem}{Theorem}

\newtheorem {lemma}[theorem]{Lemma}

\newtheorem {proposition}[theorem]{Proposition}

\newenvironment {proof}[1][Proof]{\noindent \textbf {#1.} }{\ \rule {0.5em}{0.5em}}
\begin{document}
\title{Defining Equations of Nilpotent Orbits for Borel Subgroups of Modality Zero in Type A\textsubscript {n}}

\author{
Madeleine Burkhart and David Vella}
\date{
August 15, 2017}
\maketitle

\begin{abstract}Let
$G$
be a quasi-simple algebraic group defined over an algebraically closed field
$k$
and
$B$
a Borel subgroup of
$G$
acting on the nilradical
$\mathfrak{n}$
of its Lie algebra
$\mathfrak{b}$
via the Adjoint representation.  It is known that
$B$
has only finitely many orbits in only five cases: when
$G$
is of type
$A_{n}$
for
$n \leq 4 ,$
and
$G$
is type
$B_{2} .$
In this paper, we elaborate on this work in the case when
$G =SL_{n +1}(k)$
(type $A_{n}) ,$ for $n \leq 4 ,$ by finding the defining equations of each orbit.  Consequences of these equations include the dimension of the orbits and the closure ordering on the set of orbits, although these facts are already known. The other case, when
$G$
is type $B_{2} ,$ can be approached the same way and is treated in a separate paper, where we believe the determination of the closure order is new.

\end{abstract}

\section{
Introduction.}
Let $k$ be an algebraically closed field, and $G$ a quasi-simple algebraic group defined over $k .$
Fix a maximal torus
$T$
of
$G$
, and let
$\Phi $
denote the root system of
$G$
relative to
$T$
($\Phi $
is irreducible since
$G$
is quasi-simple.)
Fix a set $\Delta $
of simple roots in
$\Phi  ,$
with corresponding set of positive roots
$\Phi ^{ +} ,$
and let
$B =TU$
(
$U$
is the unipotent radical of
$B)$
be th--SW----SW--e Borel subgroup of
$G$
determined by
$\Phi ^{ +} .$
Write the one-dimensional unipotent root group corresponding to a root
$\alpha $
as
$U_{\alpha } .$
Denote the Lie algebra of
$G$
by
$\mathfrak{g}$
, that of
$T$
by
$\mathfrak{h}$
, and that of
$B$
by
$\mathfrak{b} .$
Then the nilradical
$\mathfrak{n} =\mathfrak{n}(\mathfrak{b})$
of
$\mathfrak{b}$
is also the Lie algebra of
$U$, and we have decompositions
$\mathfrak{b} =\mathfrak{h} \oplus \mathfrak{n} ,$
and
$\mathfrak{n} = \oplus _{\alpha  \in \Phi ^{ +}}\mathfrak{g}_{\alpha }$
as vector spaces, where
$\mathfrak{g}_{\alpha }$
is the root space of
$\mathfrak{g}$
corresponding to
$\alpha  ,$
and is also the Lie algebra of
$U_{\alpha } .$.

$G$
acts on
$\mathfrak{g}$
via the Adjoint representation, and the study of the orbits of this action is a classical part of Lie theory.  It is known that there are only finitely many nilpotent
$G$-orbits (a \emph{nilpotent orbit} means an orbit of a nilpotent element of
$\mathfrak{g} .)$
There are combinatorial indexing sets for these nilpotent orbits, and there are formulas to
compute the dimension of each orbit.  Also, it is known which orbits are in the Zariski closures of
any given orbit (the \emph{closure ordering}.) Therefore, it is well understood
how all the nilpotent orbits fit together to form a larger object, called the
\emph{nullcone}
$\mathcal{N}$
of
$\mathfrak{g} ,$
which is the union of the nilpotent orbits.  Details of this classical theory can be found in [CM]
for the characteristic zero case and [C] or [J] more generally.

Beginning in the late 1980's, there has been an interest in generalizing this study to the nilpotent orbits of certain subgroups $H$ of $G ,$ particularly when $H$ is a Borel subgroup $B$ of $G$ (see [BH] and [K]), and more generally, when $H$ is a parabolic subgroup $P$ (see [BHRZ], [HR], [P], [PR], [R1], and [R2], for example.)  In [K], the following theorem is proved in case the characteristic of $k$ is zero:

\begin{theorem}
(Kashin [K], 1990)  \label{Kashin}Let $G$ be a quasi-simple over $k ,$ where $char(k) =0 ,$ and suppose $B$ is a Borel subgroup of $G$ acting on the nilradical $\mathfrak{n}$ of the Lie algebra $\mathfrak{b}$ of $B$ via the adjoint representation.  The number of orbits of $B$ on $\mathfrak{n}$ is finite (i.e., $B$ has ``modality zero'') if and only if
$G$ is type $A_{n}$ for $n \leq 4$ or $G$ is type $B_{2} .$
\end{theorem}

Over the next decade, this result was gradually extended in various ways to parabolic subgroups, and ultimately, Theorem 1.1 of [HR] related the status of a parabolic subgroup of $G$ having a finite number of orbits on the nilradical of it's Lie algebra to the nilpotency class of the unipotent radical of $P .$  From this theorem one can recover Kashin's result (Theorem \vref{Kashin}) with the added bonus that the proof is valid if $char(k) =p$ is a good prime for $G ,$ as well as characteristic zero.

Once one knows that the number of orbits is finite, one would like to mimic what is known for nilpotent $G$-orbits; that is, classify them using some index set, compute their dimensions, and determine the closure order.  In the five cases given by Kashin's theorem where $B$ has modality zero, the groups are all classical matrix groups of low rank, so we have been able to use elementary matrix calculations to do this.  Our approach is not especially elegant, but it is very detailed - we can give the defining polynomial equations for each orbit, exhibiting the orbit explicitly as an intersection of an open set and a closed set in $\mathfrak{n} .$  From there, it is easy to determine the dimensions of the orbits and the closure order.

At the time in 2015 when we worked out these defining equations, we were unaware of the paper [BHRZ], which used very different techniques in the four cases of Borel subgroups of type $A_{n}$ which have modality zero.  The closure order in these cases (as well as in some cases for parabolic subgroups) is treated there.  Since that paper treated type $A$ only, the lone case of $G =SO_{5}(k)$ (type $B_{2}$) was the last remaining case yet to appear in print.

In [BV], we have therefore written up the case for type $B_{2} ,$ where we believe the determination of the closure order and the dimensions of the orbits is new (although the closure order in this case was known to the second author as far back as 2004, it was never published.)  We also include in that paper some motivating remarks to explain why one might be interested in nilpotent $B$-orbits.

This paper contains the remainder of our work, the determination of the defining equations for the nilpotent $B$-orbits in the four cases of type $A_{n}$ when we have modality zero, as well as the dimensions and closure orders.  Because of the overlap with the conclusions in [BHRZ], it seemed less urgent to publish this work than for the case of $B_{2} .$  On the other hand, it seems at least plausible that in some application, knowing the actual defining equations of the orbits may be useful, and as far as we know, this is the only place where these equations are recorded.  For now, these results will be uploaded to the archive at http://arxiv.org, so the tables we produced are publicly available. \bigskip 

\section{
Nilpotent $B$-Orbits in Type $_{}A_{n} .$}
The results in this section for type $A$ are valid without any assumption on the characteristic of the field $k$ (which is consistent with the fact that all primes are good for type $A) .$

Let
$f$
be a polynomial in the coordinate ring of the affine space $\mathfrak{n}$. The zero set of
$f$
is written as
$Z(f)$
and
$Z(f ,g) =Z(f) \cap Z(g)$
is the set of common zeros of polynomials
$f$
and
$g .$
Then if we have a finite set of polynomials then
$Z(f_{1} ,f_{2} , . . . ,f_{k})$
is a Zariski-closed set, that is, it is an affine variety contained in $\mathfrak{n}$.  The notation
$V(f)$
denotes the complement of
$Z(f)$
- the set of non-zeros of
$f ,$
and so
$V(f)$
is a Zariski open set.  A locally closed set is an intersection of an open set and a closed
set, and in this section the orbits will turn out to be locally closed sets of the form
$V =Z(f_{1} ,f_{2} , . . . ,f_{k}) \cap V(g_{1}) \cap V(g_{2}) \cap  . . . \cap V(g_{\ell })$
for polynomials
$f_{i}$
and
$g_{j} .$
Observe that the closure of
$V$
is then
$Z(f_{1} ,f_{2} , . . . ,f_{k})$
and
$V$
is open and dense in this closure, whence
$\dim V =\dim Z(f_{1} ,f_{2} , . . . ,f_{k}) =\dim \mathfrak{n} -k ,$ provided that the $f_{i}$ are algebraically independent.  To save space, we will also abbreviate the intersection $V(g_{1}) \cap V(g_{2}) \cap  . . . \cap V(g_{\ell })$ of
open sets by $V(g_{1} ,g_{2} , . . ,g_{\ell })$.

If
$U\gamma $
is a root group of
$G ,$
then
$U_{\gamma }(t)$
denotes the image of
$t$
under the standard isomorphism
$k_{add} \approx U_{\gamma }$
. In classical groups, the Adjoint action on the Lie algebra is simply conjugation of
matrices. The matrix
$e_{ij}$
is the matrix with a
$1$
in the
$ij$
position and
$0$
everywhere else. In what follows, it
will be helpful to remember how unipotent groups act on weight vectors in rational
$G$
-modules:

\begin{lemma}
 ([H2], Proposition 27.2) Let
$\alpha  \in \Phi  ,$
and let
$v \in V_{\lambda }$
be a weight vector in any rational
$G$-module. Then each element
$u \in U_{\alpha }$
acts on
$v$
as follows:
$u .v =v +\sum _{k >0}v_{\lambda  +k\alpha } ,$
where
$v_{\lambda  +k\alpha }$
is a weight vector of weight
$\lambda  +k\alpha  ,$
and
$k$
is a positive integer.
\label{unipotentaction}
\end{lemma}

Finally, we observe that since isogenous groups all have the same orbits, we may as well assume $G =SL_{n +1}(k)$ for type $A_{n} .$

\subsection{
Type $A_{1} .$}
Here, $G =SL_{2}(k) ,$ and $\mathfrak{g} =\mathfrak{s}\mathfrak{l}_{2}(k)$.  The root system is $\Phi  =\{ \pm \alpha \} ,$ and $\Delta  =\Phi ^{ +} =\{\alpha \} .$  Thus $\mathfrak{n} =\mathfrak{g}_{\alpha }$ is one dimensional, isomorphic to the affine line as a $k$-variety.  We take the root vector to be $x =x_{\alpha } =e_{12} ,$ with corresponding coordinate function denoted $X_{\alpha } ,$ or simply $X .$  Thus, given an arbitrary element \thinspace $zx =\left [\begin{array}{cc}0 & z \\
0 & 0\end{array}\right ]$$ \in \mathfrak{g}_{\alpha }$, the function $X$ justs picks out the coordinate:\ $X(zx) =z .$  A general element of $B$ has the form $\left [\begin{array}{cc}p & b \\
0 & p^{ -1}\end{array}\right ]$ for $b ,p \in k$ and $p \neq 0.$  By direct calculation we have\begin{equation*}\left [\begin{array}{cc}p & b \\
0 & p^{ -1}\end{array}\right ]\left [\begin{array}{cc}0 & z \\
0 & 0\end{array}\right ]\left [\begin{array}{cc}p & b \\
0 & p^{ -1}\end{array}\right ]^{ -1} =\left [\begin{array}{cc}0 & zp^{2} \\
0 & 0\end{array}\right ]
\end{equation*}

Thus, if $z =0 ,$ then regardless of $p$ or $b ,$ the right side is also $0 ,$ whence $\{0\}$ is, of course, in an orbit by itself.  On the other hand, if $z \neq 0 ,$ then since $k$ is algebraically closed, then $\sqrt{z}$ exists and is nonzero in $k .$  Then the above calculation implies:\begin{equation*}\left [\begin{array}{cc}\sqrt{z} & 0 \\
0 & \sqrt{z}^{ -1}\end{array}\right ]\left [\begin{array}{cc}0 & 1 \\
0 & 0\end{array}\right ]\left [\begin{array}{cc}\sqrt{z} & 0 \\
0 & \sqrt{z}^{ -1}\end{array}\right ]^{ -1} =\left [\begin{array}{cc}0 & ^{}z \\
0 & 0\end{array}\right ]
\end{equation*}

which shows that every nonzero element of $\mathfrak{g}_{\alpha }$ is in one orbit, the orbit of the root vector $x =x_{\alpha } .$  That is, we have shown the $B .\{0\} =\{0\}$ and $B .x_{\alpha } =\mathfrak{n} -\{0\}$.  In particular, $B .\{0\}$ has dimension $0 ,$ and $B .x_{\alpha }$ is open and dense in $\mathfrak{n}$, so has dimension $1$.  Also, $\overline{B .x_{\alpha }} =\mathfrak{n} ,$ which contains $B .\{0\} ,$ so the closure order is determined.  The following Proposition summarizes these results:

\begin{proposition}
Let $G$ have type $A_{1}$ over an algebraically closed field of arbitrary characteristic, and let $B$ be a Borel subgroup of $G .$  Then $B$ has two orbits on $\mathfrak{n} ,$ the nilradical of $\mathfrak{b} .$ In particular, $B$ has modality $0.$  The table and figure below give the defining equations of the orbits, their dimensions, and the closure ordering. \label{A1orbits}

\begin{center}
\begin{tabular}[c]{|l|l|l|l|}\hline
$x$ & Defining Equations for $B .x$ & $\overline{B .x}$ & Dimension of $B .x$ \\
\hline
\hline
$0$ & $Z(X)$ & $Z(X)$ & $0$ \\
\hline
$x_{\alpha }$ & $V(X)$ & $\mathfrak{n}$ & $1$ \\
\hline
\end{tabular}\end{center}\par
\end{proposition}

\begin{center}\includegraphics[ width=0.9479166666666667in,]{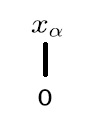}
\end{center}\par
\begin{center}
The Hasse diagram for the closure order on \end{center}\par
\begin{center}nilpotent $B$-orbits in type $A_{1}$\end{center}\par

In the Hasse diagram of the closure order, each orbit is represented by its canonical element $x$ from the first column of the table above.  \bigskip 

\subsection{
Type $A_{2} .$}
Now let $G =SL_{3}(k)$ and $\mathfrak{g} =\mathfrak{s}\mathfrak{l}_{3}(k)$.  Here $\Delta  =\{\alpha _{1} ,\alpha _{2}\} ,$ $\Phi ^{ +} =\{\alpha _{1} ,\alpha _{2} ,\alpha _{1} +\alpha _{2}\}$ and $\mathfrak{n} \approx k^{3}$.  The three root vectors are $x_{\alpha _{1}} =e_{12} ,$ $x_{\alpha _{2}} =e_{23} ,$ and $x_{\alpha _{1} +\alpha _{2}} =e_{13} .$  To make the notation less cluttered with subscripts, we abbreviate $\alpha _{1} +\alpha _{2}$ by $\alpha _{12} ,$ and more generally, in type $A_{n} ,$ every positive root has the form $\alpha _{i} +\alpha _{i +1} +\alpha _{i +2} + . . . +\alpha _{j}$ for $1 \leq i \leq j \leq n ,$ which we will abbreviate by $\alpha _{ij} .$  We use similar abbreviations for the root vectors in the corresponding root spaces, so the root vector for $\alpha _{ij}$ is just denoted $x_{ij}$ instead of $x_{\alpha _{ij}} ,$ and similarly for the coordinate functions in the coordinate ring of $\mathfrak{n}$, so the coordinate function of $x_{ij}$ is denoted $X_{ij}$ instead of $X_{\alpha _{ij}} .$  Thus, for type $A_{2} ,$ $\Phi ^{ +} =\{\alpha _{1} ,\alpha _{2} ,\alpha _{12}\} ,$ with corresponding root vectors $x_{1} ,$ $x_{2} ,$ and $x_{12} ,$ and coordinate functions $X_{1} ,X_{2} ,$ and $X_{12} .$  

We take the maximal torus $T$ to consist of diagonal matrices, so a typical element of $T$ looks like\begin{equation*}T(r ,s) =diag(r ,s ,(rs)^{ -1}) =\left [\begin{array}{ccc}r & 0 & 0 \\
0 & s & 0 \\
0 & 0 & (rs)^{ -1}\end{array}\right ]
\end{equation*}

for $r ,s \neq 0$ in $k .$  Finally, in type $A ,$ the exponential map $\exp  :g \rightarrow G$ takes $e_{ij}$ ($i \neq j)$ to $I +e_{ij}$ for each $i ,j .$  Thus the canonical isomorphism of $k_{add}$ to the root space $U_{\alpha _{ij}}$ has the property that $U_{\alpha }(r) =\exp (rx_{ij}) =I +rx_{ij} .$  For example, a typical element of $U_{\alpha _{12}}$ looks like\begin{equation*}U_{a_{12}}(r) =I +re_{13} =\left [\begin{array}{ccc}1 & 0 & r \\
0 & 1 & 0 \\
0 & 0 & 1\end{array}\right ]
\end{equation*}

In any type, $0$ is in an orbit by itself (the only closed orbit.), so $B.0 =Z(X_{1} ,X_{2} ,X_{12}) =\{0\} .$  Next, consider the orbit of the root vector $x_{12}$ corresponding to the high root $\alpha _{12}$.  By lemma \vref{unipotentaction}, every root space $U_{\alpha }$ for $\alpha  \in \Phi ^{ +}$ fixes $x_{12} .$  (Or, just multiply it out and see, or, just recall that the high root vector is a maximal vector in the Adjoint representation).  Thus, $B .x_{12} =T .x_{12} \subseteq \mathfrak{g}_{\alpha _{12}}$.  If $z \neq 0$ in $k ,$ let $d =k^{\frac{1}{3}} ,$ the cube root of $k ,$ which lies in $k$ since we assume $k$ is algebraically closed.  Then:\label{A2highroot}
\begin{gather}\left [\begin{array}{ccc}d & 0 & 0 \\
0 & d & 0 \\
0 & 0 & d^{ -2}\end{array}\right ]\left [\begin{array}{ccc}0 & 0 & 1 \\
0 & 0 & 0 \\
0 & 0 & 0\end{array}\right ]\left [\begin{array}{ccc}d & 0 & 0 \\
0 & d & 0 \\
0 & 0 & d^{ -2}\end{array}\right ]^{ -1} \nonumber  \\
 =\left [\begin{array}{ccc}0 & 0 & d^{3} \\
0 & 0 & 0 \\
0 & 0 & 0\end{array}\right ] =\left [\begin{array}{ccc}0 & 0 & z \\
0 & 0 & 0 \\
0 & 0 & 0\end{array}\right ] \label{A2highroot}\end{gather}

This shows any nonzero element of $\mathfrak{g}_{\alpha _{12}}$ is part of the orbit of $x_{12} .$  Thus, $B .x_{12} =Z(X_{1} ,X_{2}) \cap V(x_{12}) .$  Like the $A_{1}$ example, this orbit is the minimal nonzero one, and has dimension $1.$  Indeed, this is always the case for the orbit of a maximal vector in any type,

Next, consider the orbit of the root vector $x_{1} .$  By lemma \vref{unipotentaction}, every root group $U_{\alpha }$ for $\alpha $ positive fixes $x_{1}$ except for $U_{2} =U_{\alpha _{2}} .$  Thus, $B .x_{1} =TU .x_{1} =TU_{2} .x_{1}$  Since $T$ normalizes $U$ and $U_{2} ,$ we can compute these orbits one step at a time.  A direct calculation gives the $U$-orbit:
\begin{gather}U_{2}(t) .x_{1} =\left [\begin{array}{ccc}1 & 0 & 0 \\
0 & 1 & t \\
0 & 0 & 1\end{array}\right ]\left [\begin{array}{ccc}0 & 1 & 0 \\
0 & 0 & 0 \\
0 & 0 & 0\end{array}\right ]\left [\begin{array}{ccc}1 & 0 & 0 \\
0 & 1 &  -t \\
0 & 0 & 1\end{array}\right ] \nonumber  \\
 =\left [\begin{array}{ccc}0 & 1 &  -t \\
0 & 0 & 0 \\
0 & 0 & 0\end{array}\right ] =x_{1} -tx_{12} \nonumber \end{gather}

If one prefers, instead of writing elements as linear combinations of the root vectors, one could instead exploit the isomorphism of $\mathfrak{n}$ with $k^{3}$, and use the coordinate form, so $ax_{1} +bx_{2} +cx_{12} =(a ,b ,c) ,$ the $(X_{1} ,X_{2} ,X_{12})$-coordinates.  So the result of the last calculation can be expressed as $U_{2}(t) .(1 ,0 ,0) =(1 ,0 , -t)$.  Now, since orbits are conical varieties, we would expect to obtain all nonzero multiples of $(1 ,0 , -t) ,$ for any $t \in k$ in the orbit.  In other words, we would expect to obtain all elements of the form $(a ,0 ,c)$ where $a \neq 0.$ within the orbit.  That is, we would expect $Z(X_{2}) \cap V(X_{1}) \subseteq B .x_{1} .$  In fact, we can verify that this set is equal to the orbit. Containment one way is obtained by letting $T$ act on $U_{2}(t) :$\label{A2alpha1orbit}

\begin{gather}T(r ,s) .U_{2}(t) .x_{1} =T(r ,s) .(x_{1} -tx_{12}) \nonumber  \\
 =\left [\begin{array}{ccc}r & 0 & 0 \\
0 & s & 0 \\
0 & 0 & (rs)^{ -1}\end{array}\right ]\left [\begin{array}{ccc}0 & 1 &  -t \\
0 & 0 & 0 \\
0 & 0 & 0\end{array}\right ]\left [\begin{array}{ccc}r & 0 & 0 \\
0 & s & 0 \\
0 & 0 & (rs)^{ -1}\end{array}\right ]^{ -1} \nonumber  \\
 =\left [\begin{array}{ccc}0 & rs^{ -1} & r^{2}st \\
0 & 0 & 0 \\
0 & 0 & 0\end{array}\right ] \label{A2alpha1orbit}\end{gather}

Since $rs^{ -1} \neq 0 ,$ this shows $B .x_{1} \subseteq Z(X_{2}) \cap V(X_{1}) .$  Conversely, if $xx_{1} +zx_{12} =(x ,0 ,z)$ is an arbitrary element of $Z(X_{2}) \cap V(X_{1}) ,$ (so $x \neq 0) ,$ then by (\vref{A2alpha1orbit}) observe that the following element of $B$ (actually, of $TU_{2})$ sends $x_{1}$ to this element:
\begin{gather}\left [\begin{array}{ccc}x & 0 & 0 \\
0 & 1 &  -zx^{ -2} \\
0 & 0 & x^{ -1}\end{array}\right ]\left [\begin{array}{ccc}0 & 1 & 0 \\
0 & 0 & 0 \\
0 & 0 & 0\end{array}\right ]\left [\begin{array}{ccc}x & 0 & 0 \\
0 & 1 &  -zx^{2} \\
0 & 0 & x^{ -1}\end{array}\right ]^{ -1} \nonumber  \\
 =\left [\begin{array}{ccc}0 & x & z \\
0 & 0 & 0 \\
0 & 0 & 0\end{array}\right ] \nonumber \end{gather}

so $(x ,0 ,z)$ belongs to the orbit.  This shows the reverse containment, so \thinspace $B .x_{1} =Z(X_{2}) \cap V(X_{1}) .$  Note that its closure is $Z(X_{2})$ which is dimension $2$ (it's the sum of the two root spaces $\mathfrak{g}_{\alpha _{1}} \oplus \mathfrak{g}_{\alpha _{12}}$), whence $B .x_{1}$ is $2$ dimensional also.

By symmetry, or by similar calculations, we also have $B .x_{2} =Z(X_{1}) \cap V(X_{2}) ,$ which is also a $2$ dimensional orbit, dense in the sum $\mathfrak{g}\ensuremath{\operatorname*{}}_{\alpha _{2}} \oplus \mathfrak{g}_{\alpha _{12}}$.

We now consider the orbit of the sum of the root vectors for the simple roots, $B ,(x_{1} +x_{2}) .$  By lemma \vref{unipotentaction} the $U$-orbit of $x_{1} +x_{2}$ is all elements of the form $x_{1} +x_{2} +tx_{12.}$ (the reader should check this.)  Now let $T$ act to complete the $B$-orbit:

\begin{equation}T(r ,s) .(x_{1} +x_{2} +tx_{12}) =rs^{ -1}x_{1} +rs^{2}x_{2} +r^{2}stx_{12}
\end{equation}

and since $r$ and $s$ are nonzero, this shows $B .(x_{1} +x_{2}) \subseteq V(X_{1}) \cap V(X_{2}) =V(X_{1} ,X_{2}) .$  We claim this containment is actually an equality.  Indeed, let $x ,y ,z$ be elements of $k$ with $x ,y \neq 0$, so that $(x ,y ,z)$ is an arbitrary element of $V(X_{1} ,X_{2}) .$  Then observe that:
\begin{gather}\left [\begin{array}{ccc}x^{\frac{2}{3}}y^{\frac{1}{3}} & zx^{ -\frac{1}{3}}y^{ -\frac{2}{3}} & 0 \\
0 & x^{ -\frac{1}{3}}y^{\frac{1}{3}} & 0 \\
0 & 0 & x^{ -\frac{1}{3}}y^{ -\frac{2}{3}}\end{array}\right ]\left [\begin{array}{ccc}0 & 1 & 0 \\
0 & 0 & 1 \\
0 & 0 & 0\end{array}\right ]\left [\begin{array}{ccc}x^{\frac{2}{3}}y^{\frac{1}{3}} & zx^{ -\frac{1}{3}}y^{ -\frac{2}{3}} & 0 \\
0 & x^{ -\frac{1}{3}}y^{\frac{1}{3}} & 0 \\
0 & 0 & x^{ -\frac{1}{3}}y^{ -\frac{2}{3}}\end{array}\right ]^{ -1} \nonumber  \\
 =\left [\begin{array}{ccc}0 & x & z \\
0 & 0 & y \\
0 & 0 & 0\end{array}\right ] \nonumber \end{gather}

showing that an arbitrary element of $V(X_{1}) \cap V(X_{2})$ is contained in the orbit.  We have shown that $B .(x_{1} +x_{2}) =V(X_{1}) \cap V(X_{2}) =V(X_{1} ,X_{2}) .$  Notice this orbit is open and dense in $\mathfrak{n} ,$ so is $3$ dimensional, the largest possible orbit, usually called the \emph{regular} orbit.

There are no further orbits to consider.  Indeed, take an arbitrary element $(x ,y ,z)$ of $\mathfrak{n}$.  If both $x$ and $y$ are nonzero, then this element belongs to $V(X_{1} ,X_{2}) =B .(x_{1} +x_{2}) .$  If $x \neq 0$ but $y =0 ,$ then this element belongs to $Z(X_{2}) \cap V(X_{1}) =B .x_{1} ,$ while if $x =0$ and $y \neq 0 ,$ then it belongs to $Z(X_{1}) \cap V(X_{2}) =B .x_{2} .$ If both $x$ and $y$ are $0 ,$ then this element $(0 ,0 ,z)$ belongs to $B .x_{12}$ if $z \neq 0 ,$ and $B.0$ if $z =0.$   This exhausts the possibilities for $\mathfrak{n} .$

Summarizing, we have proved:

\begin{proposition}
Let $G$ have type $A_{2}$ over an algebraically closed field of arbitrary characteristic, and let $B$ be a Borel subgroup of $G .$  Then $B$ has five orbits on $\mathfrak{n} ,$ the nilradical of $\mathfrak{b} .$ In particular, $B$ has modality $0.$  The table and figure below give the defining equations of the orbits, their dimensions, and the closure ordering. \label{A2orbits}

\begin{center}
\begin{tabular}[c]{|l|l|l|l|}\hline
$x$ & Defining Equations for  $B .x$ & $\overline{B .x}$ & Dimension of \thinspace $B .x$ \\
\hline
\hline
$0$ & $Z(X_{1} ,X_{2 ,}X_{12})$ & $Z(X_{1} ,X_{2 ,}X_{12})$
& $0$ \\
\hline
$x_{12}$ & $Z(X_{1} ,X_{2}) \cap V(X_{12})$ & $Z(X_{1} ,X_{2})$ & $1$ \\
\hline
$x_{1}$ & $Z(X_{2}) \cap V(X_{1})$ & $Z(X_{2})$ & $2$ \\
\hline
$x_{2}$ & $Z(X_{1}) \cap V(X_{2})$ & $Z(X_{1})$ & $2$ \\
\hline
$x_{1} +x_{2}$ & $V(X_{1} ,X_{2})$ & $\mathfrak{n}$ & $3$ \\
\hline
\end{tabular}\end{center}\par
\end{proposition}

\begin{center}\includegraphics[ width=1.7083333333333335in,]{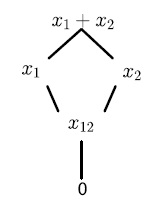}
\end{center}\par
\begin{center}The Hasse diagram for the closure order on \end{center}\par
\begin{center}nilpotent $B$-orbits in type $A_{2}$\end{center}\par

\begin{proof}
We have already verified all the entries in the table, including the dimensions, and we have also already verified that these $5$ orbits exhaust $\mathfrak{n} .$  It only remains to verify the
Hasse diagram for the closure order.  But we can see this follows from the third column in the table, since $Z(X_{1} ,X_{2} ,X_{12}) \subseteq Z(X_{1} ,X_{2}) \subseteq Z(X_{i}) \subseteq \mathfrak{n}$ for $i =1 ,2.$  This completes the proof. \bigskip 
\end{proof}

\subsection{
Type $A_{3} .$}
Now suppose $G =SL_{4}(k) .$  In this case, we shall see that there are $16$ orbits, and in the last case (type $A_{4})$ there will be $61$ orbits.  To save space, we will suppress many of the details in cases where it is straightforward to check the results, employing the same techniques we have used in the first two cases.

Here, $\Delta  =\{\alpha _{1} ,\alpha _{2} ,\alpha _{3}\} ,$ and $\Phi ^{ +} =\{\alpha _{1} ,\alpha _{2} ,\alpha _{3} ,\alpha _{12} ,\alpha _{23} ,\alpha _{13}\}$, using the abbreviations outlined above for the non-simple roots.  The root vectors in $\mathfrak{n}$ are

\begin{center}
\begin{tabular}[c]{l}$x_{1} =x_{\alpha _{1}} =e_{12} ,$\hspace*{0.3in}$x_{2} =x_{\alpha _{2}} =e_{23} ,$ \hspace*{0.3in}$x_{3} =x_{\alpha _{3}} =e_{34}$ \\
$x_{12} =x_{\alpha _{12}} =x_{\alpha _{1} +\alpha _{2}} =e_{13} ,$ \hspace*{0.3in}$x_{23} =x_{\alpha _{23}} =x_{\alpha _{2} +\alpha _{3}} =e_{24}$ \\
$x_{13} =x_{\alpha _{13}} =x_{\alpha _{1} +\alpha _{2} +\alpha _{3}} =e_{14}$
\end{tabular}\end{center}\par
so $\mathfrak{n}$ is $6$ dimensional.  The coordinate functions use the same notational abbreviations, so the coordinate ring of $\mathfrak{n}$ is $k[X_{1} ,X_{2} ,X_{3} ,X_{12} ,X_{23} ,X_{13}]$, and elements of $\mathfrak{n}$ can be written as $(u ,v ,w ,x ,y ,z) \in k^{6} .$  A typical element of $T$ has the form $T(r ,s ,t) =diag(r ,s ,t ,(rst)^{ -1}) .$  For each root vector $\alpha  ,$ we have $U_{\alpha }(t) =\exp (tx_{\alpha }) =I +tx_{\alpha } .$  We state our results at the outset:

\begin{proposition}
Let $G$ have type $A_{3}$ over an algebraically closed field of arbitrary characteristic, and let $B$ be a Borel subgroup of $G .$  Then $B$ has $16$ orbits on $\mathfrak{n} ,$ the nilradical of $\mathfrak{b} .$ In particular, $B$ has modality $0.$  The table and figure below give the defining equations of the orbits and their dimensions, and the closure ordering is pictured in the Hasse diagram following the table.

\begin{center}
\begin{tabular}[c]{|c|c|c|}\hline
$x$ & Defining equations for \thinspace $B .x$ & Dimension of $B .x$ \\
\hline
\hline
$0$ & $Z(X_{1} ,X_{2} ,X_{3} ,X_{12} ,X_{23} ,X_{13})$ & $0$ \\
\hline
$x_{13}$ & $Z(X_{1} ,X_{2} ,X_{3} ,X_{12} ,X_{23}) \cap V(X_{13})$ & $1$ \\
\hline
$x_{12}$ & $Z(X_{1} ,X_{2} ,X_{3} ,X_{23}) \cap V(X_{12})$ & $2$ \\
\hline
$x_{23}$ & $Z(X_{1} ,X_{2} ,X_{3} ,X_{12}) \cap V(X_{23})$ & $2$ \\
\hline
$x_{1}$ & $Z(X_{2} ,X_{3} ,X_{23}) \cap V(X_{1})$ & $3$ \\
\hline
$x_{2}$ & $Z(X_{1} ,X_{3} ,X_{2}X_{13} -X_{12}X_{23}) \cap V(X_{2})$ & $3$ \\
\hline
$x_{3}$ & $Z(X_{1} ,X_{2} ,X_{12}) \cap V(X_{3})$ & $3$ \\
\hline
$x_{12} +x_{23}$ & $Z(X_{1} ,X_{2} ,X_{3}) \cap V(X_{12} ,X_{23})$ & $3$ \\
\hline
$x_{2} +x_{13}$ & $Z(X_{1} ,X_{3}) \cap V(X_{2} ,X_{2}X_{13} -X_{12}X_{23})$ & $4$ \\
\hline
$x_{3} +x_{12}$ & $Z(X_{1} ,X_{2}) \cap V(X_{12}X_{3} +X_{1}X_{23})$ & $4$ \\
\hline
$x_{1} +x_{23}$ & $Z(X_{2} ,X_{3}) \cap V(X_{12}X_{3} +X_{1}X_{23})$ & $4$ \\
\hline
$x_{1} +x_{3}$ & $Z(X_{2} ,X_{12}X_{3} +X_{1}X_{23}) \cap V(X_{1} ,X_{3})$ & $4$ \\
\hline
$x_{1} +x_{2}$ & $Z(X_{3}) \cap V(X_{1} ,X_{2})$ & $5$ \\
\hline
$x_{2} +x_{3}$ & $Z(X_{1}) \cap V(X_{2} ,X_{3})$ & $5$ \\
\hline
$x_{1} +x_{3} +x_{23}$ & $Z(X_{2}) \cap V(X_{1} ,X_{3} ,X_{12}X_{3} +X_{1}X_{23})$ & $5$ \\
\hline
$x_{1} +x_{2} +x_{3}$ & $V(X_{1} ,X_{2} ,X_{3})$ & $6$ \\
\hline
\end{tabular}\end{center}\par
\end{proposition}

\begin{center}\includegraphics[ width=3.135416666666667in,]{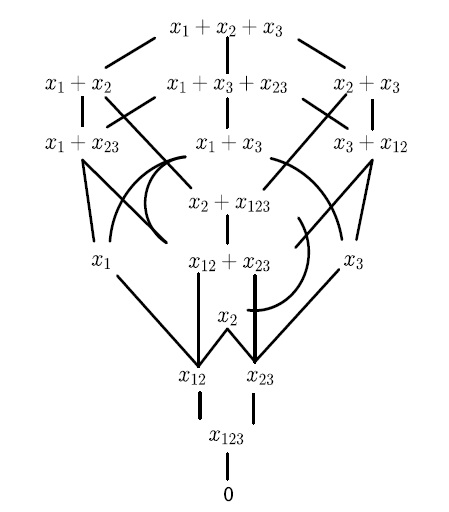}
\end{center}\par

\begin{center}The Hasse diagram for the closure order on \end{center}\par
\begin{center}
nilpotent $B$-orbits in type
$A_{3}$\textsubscript {}\end{center}\par
\begin{center}\end{center}\par

\begin{proof}
To save space, we have omitted the column containing the orbit closures from the table.  By now it is clear how to find them - they are just the closed sets in the defining equations.  Suppose for the moment that the defining equations have been established and that these $16$ orbits have been shown to be all of them.  Then the orbit closures are easy to find as remarked, and the containments of the orbit closures which give the closure ordering are easy to check, even considering the $6$ orbits with nonlinear defining polynomials. We illustrate with one example and leave the others for the reader.  

Consider \thinspace $\overline{B .(x_{12} +x_{23})} =Z(X_{1} ,X_{2} ,X_{3}) .$  Since $X_{1} =0$ and $X_{3} =0$ for points in this orbit, it follows that $X_{12}X_{3} +X_{1}X_{23} =0$ as well, whence\begin{equation*}B .(x_{12} +x_{23}) \subseteq \overline{B(x_{12} +x_{23})} \subseteq Z(X_{2} ,X_{12}X_{3} +X_{1}X_{23}) =\overline{B .(x_{1} +x_{3})} .
\end{equation*}

Therefore, we draw an edge from $B(x_{12} +x_{23})$ to $B .(x_{1} +x_{3})$ in the Hasse diagram.  All the other edges in the Hasse diagram follow similarly.

It remains to verify the defining equations and show these orbits exhaust the possibilities.  As in the above cases, we can always show that the orbit of the element in the first column is contained within the locally closed set in the same row of the second column by a simple matrix calculation, using the semidirect product decomposition $B =TU$ and lemma \vref{unipotentaction} to simplify the calculations.  For example, consider $B .x_{2}$.  By lemma \vref{unipotentaction} all root groups $U_{\gamma }$ fix $x_{2}$ except when $\gamma  =\alpha _{1}$ or $\alpha _{3} .$  Thus, we need to check the action under $U_{1}$ and $U_{3}$ (which commute with each other), and finally under $TU_{1}U_{3} .$  We have by direct matrix calculation:\label{A3alpha2orbit}\begin{equation}T(r ,s ,t)U_{1}(a)U_{3}(c) .x_{2} =st^{ -1}x_{2} +art^{ -1}x_{12} -crs^{2}tx_{23} -acr^{2}stx_{13} \label{A3alpha2orbit}
\end{equation}

But observe for the element on the right side of the above, we have $X_{2}X_{13} -X_{12}X_{23} =(st^{ -1})( -acr^{2}st) -(art^{ -1})( -crs^{2}t) =0 ,$ and since $st^{ -1} \neq 0 ,$ this yields that $B .x_{2} \subseteq Z(X_{1} ,X_{3} ,X_{2}X_{13} -X_{12}X_{23}) \cap V(X_{2}) .$  In each of the $16$ orbits, a similar calculation will give the orbit $B .x$ is contained in the locally closed set $S$ given in the second column.

The reverse containment is less trivial to demonstrate.  In each case we must take an arbitrary element $m$ of the locally closed set $S$ and  show that $m =b .x$ for some $b \in B .$  To determine such an element, we must find the solution to a nonlinear polynomial system of equations (note that we do not need the complete solution - one solution is sufficient.  There may be many elements of $b$ which carry $x$ to $m .)$

We illustrate this process with the orbit at hand.  So far we have $B .x_{2} \subseteq S =Z(X_{1} ,X_{3} ,X_{2}X_{13} -X_{12}X_{23}) \cap V(X_{2})$, and we wish to show equality.
Choose an arbitrary element \thinspace $m \in S$ which has the form $vx_{2} +xx_{12} +yx_{23} +zx_{13} =(0 ,v ,0 ,x ,y ,z) ,$ where $v \neq 0$ and $vz -xy =0.$  By (\vref{A3alpha2orbit}) we must solve the system:

\begin{center}
\begin{tabular}[c]{l}$st^{ -1} =v$ \\
$art^{ -1} =x$ \\
$ -crs^{2}t =y$ \\
$ -acr^{2}st =z$
\end{tabular}
\end{center}\par
for $a ,c ,r ,s ,t$ in terms of $v ,x ,y ,z .$  Because we have more unknowns that equations, we expect that the solution is not unique as noted above.  There is at least one extra parameter among $a ,c ,r ,s ,t$ which is not determined by $v ,x ,y ,z ,$ and which we are free to vary as we wish.  In fact, since $vz -xy =0 ,$ there are actually $2$ extra parameters.  After some trial and error, we find a solution with $r ,t =1.$ So $s =v \neq 0 ,$ $a =x ,$ $c = -ys^{ -2} = -yv^{ -2} :$ 
\begin{gather}T(r ,s ,t)U_{1}(a)U_{3}(c) .x_{2} =T(1 ,v ,1)U_{1}(x)U_{3}( -yv^{ -2}) \nonumber  \\
 =vx_{2} +xx_{12} +yx_{23} +vxyv^{ -2}x_{13} \nonumber  \\
 =vx_{2} +xx_{12} +yx_{23} +zx_{13} \nonumber \end{gather}

The last equality since $vz -xy =0$ implies $vxyv^{ -2} =v^{2}zv^{ -2} =z .$  Thus, $B .x_{2} =S =Z(X_{1} ,X_{3} ,X_{2}X_{13} -X_{12}X_{23}) \cap V(X_{2})$
as claimed.  The other cases are similar.  We leave it to the reader to check that the orbit is contained in $S$ in each case, by a calculation similar to (\vref{A3alpha2orbit}).  To show the reverse containments, the table below gives an element $b$ of $B$ in each case which has the property that $b .x$ is an arbitrary element $m$ in $S .$  Every one of the $b$ given in the table was found by solving a nonlinear system of polynomial equations as in the above example.  In some cases, we assume the field $k$ contains square roots or fourth roots of elements, so we do use the assumption of algebraic closure of $k .$

\begin{center}
\begin{tabular}[c]{|c|c|c|}\hline
\cellcolor[HTML]{FFFFFF}$x$ & \cellcolor[HTML]{FFFFFF}$m \in S$ & \cellcolor[HTML]{FFFFFF}$b \in B$ so that $b .x =m$ \\
\hline
\hline
\cellcolor[HTML]{FFFFFF}$0$ & \cellcolor[HTML]{FFFFFF}$(0 ,0 ,0 ,0 ,0 ,0)$ & \cellcolor[HTML]{FFFFFF}any $b \in B$ \\
\hline
\cellcolor[HTML]{FFFFFF}$x_{13}$ & \cellcolor[HTML]{FFFFFF}
{\begin{tabular}[c]{|l|}\hline
$(0 ,0 ,0 ,0 ,0 ,z)$ \\
\hline
$z \neq 0$
\\
\hline
\end{tabular}}
& \cellcolor[HTML]{FFFFFF}$T(1 ,1 ,z)$ \\
\hline
\cellcolor[HTML]{FFFFFF}$x_{12}$ & \cellcolor[HTML]{FFFFFF}
{\begin{tabular}[c]{|l|}\hline
$(0 ,0 ,0 ,x ,0 ,z)$ \\
\hline
$x \neq 0$ \\
\hline
\end{tabular}} & \cellcolor[HTML]{FFFFFF}$T(x ,1 ,1)U_{3}( -zx^{ -2})$ \\
\hline
\cellcolor[HTML]{FFFFFF}$x_{23}$ & \cellcolor[HTML]{FFFFFF}
{\begin{tabular}[c]{|l|}\hline
$(0 ,0 ,0 ,0 ,y ,z)$ \\
\hline
$y \neq 0$
\\
\hline
\end{tabular}} & \cellcolor[HTML]{FFFFFF}$T(y ,1 ,1)U_{1}(zy^{ -2})$ \\
\hline
\cellcolor[HTML]{FFFFFF}$x_{1}$ & \cellcolor[HTML]{FFFFFF}
{\begin{tabular}[c]{|l|}\hline
$(u ,0 ,0 ,x ,0 ,z)$ \\
\hline
$u \neq 0$
\\
\hline
\end{tabular}} & \cellcolor[HTML]{FFFFFF}$T(u ,1 ,1)U_{23}( -zu^{ -2})U_{2}( -xu^{ -1})$ \\
\hline
\cellcolor[HTML]{FFFFFF}$x_{2}$ & \cellcolor[HTML]{FFFFFF}
{\begin{tabular}[c]{|l|}\hline
$(0 ,v ,0 ,x ,y ,z)$ \\
\hline
$v \neq 0$ \\
\hline
$vz -xy =0$ \\
\hline
\end{tabular}}\par
& \cellcolor[HTML]{FFFFFF}$T(1 ,w ,1)U_{1}(x)U_{3}( -yw^{ -2})$ \\
\hline
\cellcolor[HTML]{FFFFFF}$x_{3}$ & \cellcolor[HTML]{FFFFFF}
{\begin{tabular}[c]{|l|}\hline
$(0 ,0 ,w ,0 ,y ,z)$ \\
\hline
$w \neq 0$
\\
\hline
\end{tabular}} & \cellcolor[HTML]{FFFFFF}$T(w ,1 ,1)U_{12}(zw^{ -2})U_{2}(yw^{ -1})$ \\
\hline
\cellcolor[HTML]{FFFFFF}$x_{12} +x_{23}$ & \cellcolor[HTML]{FFFFFF}
{\begin{tabular}[c]{|l|}\hline
$(0 ,0 ,0 ,x ,y ,z)$ \\
\hline
$x ,y \neq 0$ \\
\hline
\end{tabular}} & \cellcolor[HTML]{FFFFFF}$T(x ,\sqrt{yx^{ -1}} ,1)U_{1}(yz(xy)^{ -\frac{3}{2}})$ \\
\hline
\cellcolor[HTML]{FFFFFF}$x_{2} +x_{13}$ & \cellcolor[HTML]{FFFFFF}
{\begin{tabular}[c]{|l|}\hline
$(0 ,v ,0 ,x ,y ,z)$ \\
\hline
$v \neq 0$ \\
\hline
$vz -xy \neq 0$
\\
\hline
\end{tabular}} & \cellcolor[HTML]{FFFFFF}
{\begin{tabular}[c]{|c|}\hline
$T(v^{ -1}\sqrt{vz -xy} ,v ,1) \cdot $ \\
\hline
$U_{3}( -yw^{ -1}(wz -xy)^{ -\frac{1}{2}})U_{1}(wx(wz -xy)^{ -\frac{1}{2}})$ \\
\hline
\end{tabular}} \\
\hline
\cellcolor[HTML]{FFFFFF}$x_{3} +x_{12}$ & \cellcolor[HTML]{FFFFFF}
{\begin{tabular}[c]{|l|}\hline
$(0 ,0 ,w ,x ,y ,z)$ \\
\hline
$w ,x \neq 0$
\\
\hline
\end{tabular}} & \cellcolor[HTML]{FFFFFF}$T(x ,wx^{ -1} ,1)U_{12}(z(xw)^{ -1})U_{2}(xyw^{ -2})$ \\
\hline
\cellcolor[HTML]{FFFFFF}$x_{1} +x_{23}$ & \cellcolor[HTML]{FFFFFF} 
{\begin{tabular}[c]{|l|}\hline
$(u ,0 ,0 ,x ,y ,z)$ \\
\hline
$u ,y \neq 0$
\\
\hline
\end{tabular}}
& \cellcolor[HTML]{FFFFFF}$T(u ,1 ,yu^{ -1})U_{1}(z(uy)^{ -1})U_{2}( -xyu^{ -2})$ \\
\hline
\cellcolor[HTML]{FFFFFF}$x_{1} +x_{3}$ & \cellcolor[HTML]{FFFFFF}
{\begin{tabular}[c]{|l|}\hline
$(u ,0 ,w ,x ,y ,z)$ \\
\hline
$u ,w \neq 0$ \\
\hline
$wx +uy \neq 0$ \\
\hline
\end{tabular}} & \cellcolor[HTML]{FFFFFF}
{\begin{tabular}[c]{|c|}\hline
$T(u ,1 ,\sqrt{wu^{ -1}}) \cdot $ \\
\hline
$U_{12}(zu^{ -2}\sqrt{uw^{ -1}})U_{2}(yu^{ -1}\sqrt{uw^{ -1}})$ \\
\hline
\end{tabular}} \\
\hline
\cellcolor[HTML]{FFFFFF}$x_{1} +x_{2}$ & \cellcolor[HTML]{FFFFFF} 
{\begin{tabular}[c]{|l|}\hline
$(u ,v ,0 ,x ,y ,z)$ \\
\hline
$u ,v \neq 0$ \\
\hline
\end{tabular}}
& \cellcolor[HTML]{FFFFFF}
{\begin{tabular}[c]{|c|}\hline
$T(uv ,v ,1) \cdot $ \\
\hline
$U_{23}((xy -vz)v^{ -4}u^{ -2})U_{3}( -yu^{ -1}v^{ -3})U_{1}(xv^{ -1}u^{ -1})$
\\
\hline
\end{tabular}} \\
\hline
\cellcolor[HTML]{FFFFFF}$x_{2} +x_{3}$ & \cellcolor[HTML]{FFFFFF}  
{\begin{tabular}[c]{|l|}\hline
$(0 ,v ,w ,x ,y ,z)$ \\
\hline
$v ,w \neq 0$ \\
\hline
\end{tabular}}
& \cellcolor[HTML]{FFFFFF}
{\begin{tabular}[c]{|c|}\hline
$T(v^{2}w ,1 ,v^{ -1}) \cdot $ \\
\hline
$U_{12}((vz -xy)v^{ -4}w^{ -2})U_{3}( -yv^{ -1}w^{ -1})U_{1}(xv^{ -3}w^{ -1})$ \\
\hline
\end{tabular}} \\
\hline
\cellcolor[HTML]{FFFFFF}$x_{1} +x_{3} +x_{23}$ & \cellcolor[HTML]{FFFFFF}\par   
{\begin{tabular}[c]{|l|}\hline
$(u ,0 ,w ,x ,y ,z)$ \\
\hline
$u ,w \neq 0$ \\
\hline
$wx +uy \neq 0$ \\
\hline
\end{tabular}}
& \cellcolor[HTML]{FFFFFF}
{\begin{tabular}[c]{|c|}\hline
$T((\frac{u(yuv)^{2}}{v})^{\frac{1}{4}} ,\thinspace (\frac{(yu +vx)^{2}}{u^{3}v})^{\frac{1}{4}} ,\thinspace (\frac{uv^{3}}{(yu +vx)^{2}})^{\frac{1}{4}}) \cdot $ \\
\hline
$U_{12}(\frac{vxz}{(yu +vx)^{2}})U_{1}(\frac{z}{yu +vx})U_{2}( -\frac{vx}{yu +vx})$  \\
\hline
\end{tabular}} \\
\hline
\cellcolor[HTML]{FFFFFF}$x_{1} +x_{2} +x_{3}$ & \cellcolor[HTML]{FFFFFF} 
{\begin{tabular}[c]{|l|}\hline
$(u ,v ,w ,x ,y ,z)$ \\
\hline
$u ,v ,w \neq 0$ \\
\hline
\end{tabular}}
& \cellcolor[HTML]{FFFFFF}
{\begin{tabular}[c]{|c|}\hline
$T(u^{\frac{3}{4}}v^{\frac{2}{4}}w^{\frac{1}{4}} ,\thinspace u^{ -\frac{1}{4}}v^{\frac{2}{4}}w^{\frac{1}{4}} ,\thinspace u^{ -\frac{1}{4}}v^{ -\frac{2}{4}}w^{\frac{1}{4}}) \cdot $ \\
\hline
$U_{12}((vwz -y^{2}u -wxy)u^{ -1}v^{ -2}w^{ -2}) \cdot $ \\
\hline
$U_{1}((xw +yu)(uvw)^{ -1})U_{2}(y(vw)^{ -1})$  \\
\hline
\end{tabular}} \\
\hline
\end{tabular}\end{center}\par
This will establish the reverse containments in each case and therefore yield the stated defining equations for each orbit.

Lastly, we verify that these orbits exhaust $\mathfrak{n} .$  Let $m =(u ,v ,w ,x ,y ,z)$ be an arbitrary element of $\mathfrak{n} .$We compare it to the various forms in the defining relations in the middle column of the table in the statement of the theorem.  If $u ,v ,w$ are all nonzero, then $m$ belongs to the regular orbit $B .(x_{1} +x_{2} +x_{3}) .$  So now assume at least one of $u ,v ,w$ is $0.$  If only $u =0 ,$ then $m \in Z(X_{1}) \cap V(X_{1} ,X_{2}) =B .(x_{2} +x_{3}) .$  If only $w =0 ,$ then $m \in Z(X_{3}) \cap V(X_{1} ,X_{2}) =B .(x_{1} +x_{2}) .$  If only $v =0 ,$ then there are two possibilities: $m \in B .(x_{1} +x_{3})$ if $wx +uy =0$, and $m \in B .(x_{1} +x_{3} +x_{23})$ if $wx +uy \neq 0$.

Now consider the three cases when exactly two of $u ,v ,w$ are $0.$  If $u ,v =0 ,$ there are two possibilities: $m \in B .(x_{3} +x_{12})$ if  $x \neq 0 ,$ and $m \notin B .x_{3}$ if $x =0.$ Similarly, if $v ,w =0 ,$ then $m \in B .(x_{1} +x_{23})$ if $y \neq 0$, and $m \in B .x_{1}$ if $y =0.$  If $u ,w =0 ,$ then $m \in B .x_{2}$ if $vz -xy =0 ,$ $m \in B .(x_{2} +x_{13})$ if $vz -xy \neq 0.$

That leaves the case when all three of $u ,v ,w =0.$  First, suppose $x ,y$ are both nonzero.  Then $m \in B .(x_{12} +x_{23}) ,$ so we may assume that at least one of $x$ and $y$ are $0.$  If $x =0$ and $y \neq 0 ,$ then $m \in B .x_{23} ,$ and similarly, if $x \neq 0$ and $y =0 ,$ then $m \in B .x_{12} .$  If both $x ,y =0$ then either $z \neq 0$ and $m$ is in the minimal orbit $B .x_{13} ,$ or $z =0$ and $m =0$ is in $B.0.$

This exhausts the possibilities, and completes the proof. \bigskip 
\end{proof}

\subsection{
Type $A_{4}$.}
We now assume $G =SL_{5}(k)$ and $\mathfrak{g} =\mathfrak{s}\mathfrak{l}_{5}(k) .$  Here $\Delta  =\{\alpha _{1} ,\alpha _{2} ,\alpha _{3} ,\alpha _{4}\}$ and $\Phi  =\{\alpha _{ij}$ \textbar{} $1 \leq i \leq j \leq 4\}$ (where $\alpha _{ii}$ is interpreted as just $\alpha _{i}) .$  There are $10$ positive roots, so $\mathfrak{n}$ is $10$ dimensional.  The $10$ root vectors are given by the matrices $x_{ij} =x_{\alpha _{ij}} =e_{i ,j +1} ,$ with corresponding coordinate function $X_{ij} .$  We may abbreviate an element $n =qx_{1} +rx_{2} +sx_{3} +tx_{4} +ux_{12} +vx_{23} +wx_{34} +xx_{13} +yx_{24} +zx_{14}$ of $n$ by $(q ,r ,s ,t ,u ,v ,w ,x ,y ,z) \in k^{10} .$  A typical element of $T$ has the form $T(a ,b ,c ,d) =diag(a ,b ,c ,d ,(abcd)^{ -1}) .$  For each root $\gamma  ,$ a typical element of the root group $U_{\gamma }$ has the form $U_{\gamma }(f) =$ $\exp (fx_{\gamma }) =$ $I +fx_{\gamma }$, although since $U$ is a product of $10$ such elements, we often find it more convenient to avoid lemma \vref{unipotentaction} and use the arbitrary element of $U :$
\begin{equation*}U_{arb} =\left [\begin{array}{ccccc}1 & f_{1} & f_{5} & f_{8} & f_{10} \\
0 & 1 & f_{2} & f_{6} & f_{9} \\
0 & 0 & 1 & f_{3} & f_{7} \\
0 & 0 & 0 & 1 & f_{4} \\
0 & 0 & 0 & 0 & 1\end{array}\right ]
\end{equation*}

Our results follow:

\begin{proposition}
Let $G$ have type $A_{4}$ over an algebraically closed field of arbitrary characteristic, and let $B$ be a Borel subgroup of $G .$  Then $B$ has $61$ orbits on $\mathfrak{n} ,$ the nilradical of $\mathfrak{b} .$ In particular, $B$ has modality $0.$  The tables and figure below give the defining equations of the orbits and their dimensions, and the closure ordering is pictured in the Hasse diagram following the tables. The tables organize the orbits according to their dimensions.

\begin{center}
\begin{tabular}[c]{|c|c|c|}\hline
$n$ & Defining equation of $B .n$ & Dimension of $B .n$ \\
\hline
\hline
$0$ & $Z(X_{1} ,X_{2} ,X_{3} ,X_{4} ,X_{12} ,X_{23} ,X_{34} ,X_{13} ,X_{24} ,X_{14})$ & $0$ \\
\hline
$x_{14}$ & $Z(X_{1} ,X_{2} ,X_{3} ,X_{4} ,X_{12} ,X_{23} ,X_{34} ,X_{13} ,X_{24}) \cap V(X_{14})$ & $1$ \\
\hline
$x_{13}$ & $Z(X_{1} ,X_{2} ,X_{3} ,X_{4} ,X_{12} ,X_{23} ,X_{34} ,X_{24}) \cap V(X_{13})$
& $2$ \\
\hline
$x_{24}$ & $Z(X_{1} ,X_{2} ,X_{3} ,X_{4} ,X_{12} ,X_{23} ,X_{34} ,X_{13}) \cap V(X_{24})$
& $2$ \\
\hline
\multicolumn{3}{|c|}{
Orbits of dimension at most $2$
} \\
\hline
\end{tabular}\end{center}\par
\begin{center}
\end{center}\par
\begin{center}
\begin{tabular}[c]{|c|c|}\hline
$n$ & Defining equation of $B .n$
\\
\hline
\hline
$x_{13} +x_{24}$ & $Z(X_{1} ,X_{2} ,X_{3} ,X_{4} ,X_{12} ,X_{23} ,X_{34}) \cap V(X_{13} ,X_{24})$ \\
\hline
$x_{12}$ & $Z(X_{1} ,X_{2} ,X_{3} ,X_{4} ,X_{23} ,X_{34} ,X_{24}) \cap V(X_{12})$
\\
\hline
$x_{23}$ & $Z(X_{1} ,X_{2} ,X_{3} ,X_{4} ,X_{12} ,X_{34} ,\thinspace X_{13}X_{24} -X_{23}X_{14}) \cap V(X_{23})$
\\
\hline
$x_{34}$ & $Z(X_{1} ,X_{2} ,X_{3} ,X_{4} ,X_{12} ,X_{23} ,X_{13}) \cap V(X_{34})$
\\
\hline
\multicolumn{2}{|c|}{
Orbits of dimension $3$
} \\
\hline
\end{tabular}\end{center}\par
\begin{center}
\end{center}\par
\begin{center}
\begin{tabular}[c]{|c|c|}\hline
$n$ & Defining equation of $B .n$
\\
\hline
\hline
$x_{12} +x_{24}$ & $Z(X_{1} ,X_{2} ,X_{3} ,X_{4} ,X_{23} ,X_{34}) \cap V(X_{12} ,X_{24})$
\\
\hline
$x_{23} +x_{14}$ & $Z(X_{1} ,X_{2} ,X_{3} ,X_{4} ,X_{12} ,X_{34}) \cap V(X_{23} ,\thinspace X_{13}X_{24} -X_{23}X_{14})$
\\
\hline
$x_{34} +x_{13}$ & $Z(X_{1} ,X_{2} ,X_{3} ,X_{4} ,X_{12} ,X_{23}) \cap V(X_{13} ,X_{34})$
\\
\hline
$x_{1}$ & $Z(X_{2} ,X_{3} ,X_{4} ,X_{23} ,X_{34} ,X_{24}) \cap V(X_{1})$
\\
\hline
$x_{2}$ & $Z(X_{1} ,X_{3} ,X_{4} ,X_{34} ,\thinspace X_{2}X_{13} -X_{12}X_{23} ,\thinspace X_{2}X_{14} -X_{12}X_{24}) \cap V(X_{2})$
\\
\hline
$x_{3}$ & $Z(X_{1} ,X_{2} ,X_{4} ,X_{12} ,\thinspace X_{3}X_{24} -X_{23}X_{34} ,\thinspace X_{3}X_{14} -X_{13}X_{34}) \cap V(X_{2})$
\\
\hline
$x_{4}$ & $Z(X_{1} ,X_{2} ,X_{3} ,X_{12} ,X_{23} ,X_{13}) \cap V(X_{4})$
\\
\hline
\multicolumn{2}{|c|}{
Orbits of dimension $4$
} \\
\hline
\end{tabular}\end{center}\par
\begin{center}
\begin{tabular}[c]{|c|c|}\hline
$n$ & Defining equation of $B .n$
\\
\hline
$x_{12} +x_{23}$ & $Z(X_{1} ,X_{2} ,X_{3} ,X_{4} ,X_{34}) \cap V(X_{12} ,X_{23})$
\\
\hline
$x_{12} +x_{34}$ & $Z(X_{1} ,X_{2} ,X_{3} ,X_{4} ,X_{23}) \cap V(X_{12} ,X_{34})$
\\
\hline
$x_{23} +x_{34}$ & $Z(X_{1} ,X_{2} ,X_{3} ,X_{4} ,X_{12}) \cap V(X_{23} ,X_{34})$
\\
\hline
$x_{1} +x_{24}$ & $Z(X_{2} ,X_{3} ,X_{4} ,X_{23} ,X_{34}) \cap V(X_{1} ,X_{24})$
\\
\hline
$x_{1} +x_{34}$ & $Z(X_{2} ,X_{3} ,X_{4} ,X_{23} ,\thinspace X_{12}X_{34} +X_{1}X_{24}) \cap V(X_{1} ,X_{34})$
\\
\hline
$x_{2} +x_{14}$ & $Z(X_{1} ,X_{3} ,X_{4} ,X_{34} ,\thinspace X_{12}X_{23} -X_{2}X_{13}) \cap V(X_{2} ,\thinspace X_{12}X_{24} -X_{2}X_{14})$
\\
\hline
$x_{3} +x_{14}$ & $Z(X_{1} ,X_{2} ,X_{4} ,X_{12} ,\thinspace X_{23}X_{34} -X_{3}X_{24}) \cap V(X_{3} ,\thinspace X_{13}X_{34} -X_{3}X_{14})$
\\
\hline
$x_{4} +x_{12}$ & $Z(X_{1} ,X_{2} ,X_{3} ,X_{23} ,\thinspace X_{12}X_{34} +X_{13}X_{4}) \cap V(X_{12} ,X_{4})$
\\
\hline
$x_{4} +x_{13}$ & $Z(X_{1} ,X_{2} ,X_{3} ,X_{34} ,X_{12} ,X_{23}) \cap V(X_{13} ,X_{4})$
\\
\hline
\multicolumn{2}{|c|}{Orbits of dimension $5$} \\
\hline
\end{tabular}\end{center}\par
\begin{center}
\end{center}\par
\begin{center}
\begin{tabular}[c]{|c|c|}\hline
$n$ & Defining equation of $B .n$
\\
\hline
\hline
$x_{12} +x_{23} +x_{34}$ & $Z(X_{1} ,X_{2} ,X_{3} ,X_{4}) \cap V(X_{12} ,X_{23} ,X_{34})$
\\
\hline
$x_{1} +x_{23}$ & $Z(X_{2} ,X_{3} ,X_{4} ,X_{34}) \cap V(X_{1} ,X_{23})$
\\
\hline
$x_{1} +x_{34} +x_{24}$ & $Z(X_{2} ,X_{3} ,X_{4} ,X_{23}) \cap V(X_{1} ,X_{34} ,\thinspace X_{1}X_{24} +X_{12}X_{34})$
\\
\hline
$x_{2} +x_{13}$ & $Z(X_{1} ,X_{3} ,X_{4} ,X_{34}) \cap V(X_{2} ,\thinspace X_{12}X_{23} -X_{123}X_{2})$
\\
\hline
$x_{2} +x_{34}$ & $Z(X_{1} ,X_{3} ,X_{4} ,\thinspace X_{2}X_{13} -X_{12}X_{23}) \cap V(X_{2} ,X_{34})$
\\
\hline
$x_{3} +x_{12}$ & $Z(X_{1} ,X_{2} ,X_{4} ,\thinspace X_{3}X_{24} -X_{23}X_{34}) \cap V(X_{3} ,X_{12})$
\\
\hline
$x_{3} +x_{24}$ & $Z(X_{1} ,X_{2} ,X_{4} ,X_{12}) \cap V(X_{3} ,\thinspace X_{23}X_{34} -X_{3}X_{234})$
\\
\hline
$x_{4} +x_{12} +x_{13}$ & $Z(X_{2} ,X_{3} ,X_{4} ,X_{23}) \cap V(X_{4} ,X_{12} ,\thinspace X_{4}X_{13} +X_{12}X_{34})$
\\
\hline
$x_{4} +x_{23}$ & $Z(X_{1} ,X_{2} ,X_{3} ,X_{12}) \cap V(X_{4} ,X_{23})$
\\
\hline
$x_{1} +x_{3}$ & $Z(X_{2} ,X_{4} ,\thinspace X_{1}X_{23} +X_{3}X_{12} ,X_{3}X_{24} -X_{23}X_{34}) \cap V(X_{1} ,X_{3})$
\\
\hline
$x_{1} +x_{4}$ & $Z(X_{2} ,X_{3} ,\thinspace X_{23} ,X_{1}X_{24} +X_{12}X_{34} +X_{4}X_{123}) \cap V(X_{1} ,X_{4})$
\\
\hline
$x_{2} +x_{4}$ & $Z(X_{1} ,X_{3} ,\thinspace X_{4}X_{23} +X_{2}X_{34} ,\thinspace X_{2}X_{13} -X_{12}X_{23}) \cap V(X_{2} ,X_{4})$
\\
\hline
\multicolumn{2}{|c|}{Orbits of dimension $6$
} \\
\hline
\end{tabular}\end{center}\par
\begin{center}
\end{center}\par
\begin{center}
\begin{tabular}[c]{|c|c|}\hline
$n$ & Defining equation of $B .n$
\\
\hline
\hline
$x_{1} +x_{23} +x_{34}$ & $Z(X_{2} ,X_{3} ,X_{4}) \cap V(X_{1} ,X_{23} ,X_{34})$
\\
\hline
$x_{2} +x_{34} +x_{13}$ & $Z(X_{1} ,X_{3} ,X_{4}) \cap V(X_{2} ,X_{34} ,\thinspace X_{2}X_{13} -X_{12}X_{23})$
\\
\hline
$x_{3} +x_{12} +x_{24}$ & $Z(X_{1} ,X_{2} ,X_{4}) \cap V(X_{3} ,X_{12} ,\thinspace X_{3}X_{24} -X_{23}X_{34})$
\\
\hline
$x_{4} +x_{12} +x_{23}$ & $Z(X_{1} ,X_{2} ,X_{3}) \cap V(X_{4} ,X_{12} ,X_{23})$
\\
\hline
$x_{1} +x_{2}$ & $Z(X_{3} ,X_{4} ,X_{34}) \cap V(X_{1} ,X_{2})$
\\
\hline
$x_{1} +x_{3} +x_{12}$ & $Z(X_{2} ,X_{4} ,\thinspace X_{3}X_{24} -X_{23}X_{34}) \cap V(X_{1} ,X_{3} ,\thinspace X_{1}X_{23} +X_{3}X_{12})$
\\
\hline
$x_{1} +x_{3} +x_{24}$ & $Z(X_{2} ,X_{4} ,\thinspace X_{1}X_{23} +X_{3}X_{12}) \cap V(X_{1} ,X_{3} ,\thinspace X_{3}X_{24} -X_{23}X_{34})$
\\
\hline
$x_{2} +x_{3}$ & 
{\begin{tabular}[c]{c}$Z(X_{1} ,X_{4} ,\thinspace X_{3}X_{12}X_{24} +X_{2}X_{34}X_{13} -X_{12}X_{23}X_{34} -X_{2}X_{3}X_{14})$ \\
$ \cap V(X_{2} ,X_{3} ,)$
\end{tabular}}
\\
\hline
$x_{1} +x_{4} +x_{24}$ & $Z(X_{2} ,X_{3} ,X_{23}) \cap V(X_{1} ,X_{4} ,\thinspace X_{1}X_{24} +X_{12}X_{34} +X_{4}X_{13})$
\\
\hline
$x_{2} +x_{4} +x_{13}$ & $Z(X_{1} ,X_{3} ,\thinspace X_{4}X_{23} +X_{2}X_{34}) \cap V(X_{2} ,X_{4} ,\thinspace X_{2}X_{13} -X_{12}X_{23})$
\\
\hline
$x_{2} +x_{4} +x_{34}$ & $Z(X_{1} ,X_{3} ,\thinspace X_{2}X_{13} -X_{12}X_{23}) \cap V(X_{2} ,X_{4} ,\thinspace X_{2}X_{34} +X_{4}X_{23})$
\\
\hline
$x_{3} +x_{4}$ & $Z(X_{1} ,X_{2} ,X_{12}) \cap V(X_{3} ,X_{4})$
\\
\hline
\multicolumn{2}{|c|}{Orbits of dimension $7$
} \\
\hline
\end{tabular}\end{center}\par
\begin{center}
\end{center}\par
\begin{center}
\begin{tabular}[c]{|c|c|}\hline
$n$ & Defining equation of $B .n$
\\
\hline
\hline
$x_{1} +x_{2} +x_{34}$ & $Z(X_{3} ,X_{4}) \cap V(X_{1} ,X_{2} ,X_{34})$
\\
\hline
$x_{1} +x_{3} +x_{12} +x_{24}$ & $Z(X_{2} ,X_{4}) \cap V(X_{1} ,\thinspace X_{3} ,\thinspace X_{1}X_{23} +X_{3}X_{12} ,\thinspace X_{3}X_{24} -X_{23}X_{34})$
\\
\hline
$x_{2} +x_{3} +x_{14}$ &  
{\begin{tabular}[c]{c}$Z(X_{1} ,X_{4}) \cap $ \\
$V(X_{2} ,\thinspace X_{3} ,\thinspace X_{3}X_{12}X_{24} +X_{2}X_{34}X_{13} -X_{12}X_{23}X_{34} -X_{2}X_{3}X_{14})$
\end{tabular}}
\\
\hline
$x_{1} +x_{4} +x_{23}$ & $Z(X_{2} ,X_{3}) \cap V(X_{1} ,X_{4} ,X_{23})$
\\
\hline
$x_{2} +x_{4} +x_{34} +x_{13}$ & $Z(X_{1} ,X_{3}) \cap V(X_{2} ,\thinspace X_{4} ,\thinspace X_{4}X_{23} +X_{2}X_{34} ,\thinspace X_{2}X_{13} -X_{12}X_{23})$
\\
\hline
$x_{3} +x_{4} +x_{12}$ & $Z(X_{1} ,X_{2}) \cap V(X_{3} ,X_{4} ,X_{12})$
\\
\hline
$x_{1} +x_{2} +x_{4}$ & $Z(X_{3} ,\thinspace X_{2}X_{34} +X_{4}X_{23}) \cap V(X_{1} ,X_{2} ,X_{4})$
\\
\hline
$x_{1} +x_{3} +x_{4}$ & $Z(X_{2} ,\thinspace X_{3}X_{12} +X_{1}X_{23}) \cap V(X_{1} ,X_{3} ,X_{4})$
\\
\hline
\multicolumn{2}{|c|}{Orbits of dimension $8$
} \\
\hline
\end{tabular}\end{center}\par
\begin{center}
\end{center}\par
\begin{center}
\begin{tabular}[c]{|c|c|l|}\hline
$n$ & Defining equation of $B .n$
& Dimension of $B .n$
\\
\hline
\hline
\thinspace $x_{1} +x_{2} +x_{3}$ & $Z(X_{4}) \cap V(X_{1} ,X_{2} ,X_{3})$
&
$9$
\\
\hline
$x_{1} +x_{2} +x_{4} +x_{34}$ & $Z(X_{3}) \cap V(X_{1} ,X_{2} ,X_{4} ,\thinspace X_{2}X_{34} +X_{4}X_{23})$
&
$9$
\\
\hline
$x_{1} +x_{3} +x_{4} +x_{12}$ & $Z(X_{2}) \cap V(X_{1} ,X_{3} ,X_{4} ,\thinspace X_{1}X_{23} +X_{3}X_{12})$
&
$9$
\\
\hline
$x_{2} +x_{3} +x_{4}$ & $Z(X_{1}) \cap V(X_{2} ,X_{3} ,X_{4})$
&
$9$
\\
\hline
$x_{1} +x_{2} +x_{3} +x_{4}$ & $V(X_{1} ,X_{2} ,X_{3} ,X_{4})$ &
$10$
\\
\hline
\multicolumn{3}{|c|}{Orbits of dimension $9$ or $10$
} \\
\hline
\end{tabular}\end{center}\par
\end{proposition}

\begin{center}\includegraphics[ width=5.489583333333334in,]{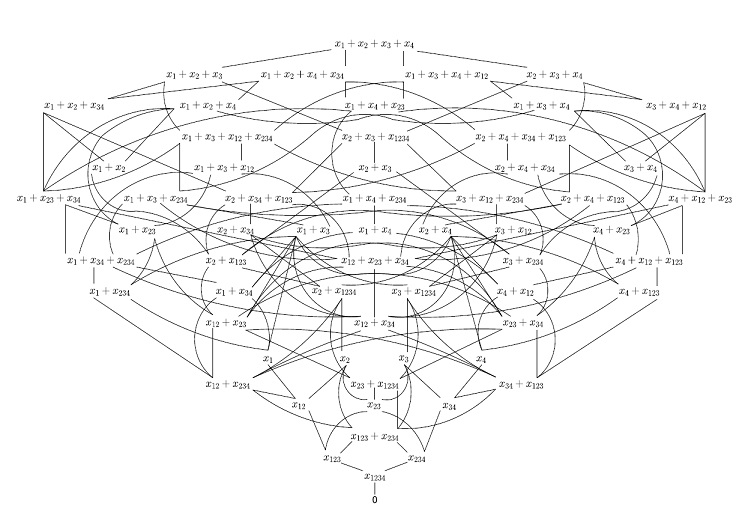}
\end{center}\par
\begin{center}The Hasse diagram for the closure order on \end{center}\par
\begin{center}nilpotent $B$-orbits in type $A_{4}$\end{center}\par

\begin{proof}
Assuming for the moment that the defining equations in the tables are correct, then note that the dimensions follow from this.  As usual, the dimension is $\dim \mathfrak{n} -$(number of algebraically independent elements of $\overline{B .x})$  We illustrate the need to be careful about the algebraic independence with an example.  Consider the orbit $B .x_{2} .$  Here we obtain for the $U$ orbit:
\begin{equation*}U_{arb} .x_{2} =x_{2} +f_{1}x_{12} -f_{3}x_{23} -f_{1}f_{3}x_{13} +(f_{3}f_{4} -f_{7})x_{24} +(f_{1}f_{3}f_{4} -f_{1}f_{7})x_{14}
\end{equation*}

For this element, $X_{2}X_{13} -X_{12}X_{23} =0 ,$ $X_{2}X_{14} -X_{12}X_{24} =0 ,$ $X_{13}X_{24} -X_{23}X_{14} =0 ,$ and $X_{2} \neq 0.$  Since all four of these conditions are homogeneous equations, and orbits are conical varieties, the entire $B$ orbit will satisfy the same equations (or just let $T(a ,b ,c ,d)$ act on this result and check the equations afterwards.)  It follows that:\label{A4alpha2OrbitWrong}\begin{equation}B .x_{2} \subseteq Z(X_{2}X_{13} -X_{12}X_{23},\thinspace X_{2}X_{14} -X_{12}X_{24} ,\thinspace X_{13}X_{24} -X_{23}X_{14}) \cap V(X_{2}) \label{A4alpha2OrbitWrong}
\end{equation}

However, the first three polynomials, which define the closure of this orbit, are algebraically dependent:\begin{equation*}X_{2}(X_{13}X_{24} -X_{23}X_{14}) =X_{24}(X_{2}X_{13} -X_{12}X_{23}) -X_{23}(X_{2}X_{14} -X_{12}X_{24})
\end{equation*}

and since $X_{2} \neq 0 ,$ it follows that if $X_{2}X_{13} -X_{12}X_{23}$ and $X_{2}X_{14} -X_{12}X_{24}$ vanish at a point in $V(X_{2}) ,$ so does $X_{13}X_{24} -X_{23}X_{14} .$  It follows that the orbit also has the following description:\label{A4alpha2OrbitRight}\begin{equation}B .x_{2} \subseteq Z(X_{2}X_{13} -X_{12}X_{23} ,\thinspace X_{2}X_{14} -X_{12}X_{24}) \cap V(X_{2}) \label{A4alpha2OrbitRight}
\end{equation}

Now the elements of the closure of this are algebraically independent, so in order to obtain the correct dimension of this closure, one must use the second description (\vref{A4alpha2OrbitRight}) which avoids the algebraic dependencies of the first one (\vref{A4alpha2OrbitWrong}).  In the tables above, the closures of the locally closed sets in the second column are all expressed in terms of algebraically independent elements of the coordinate ring of $\mathfrak{n} .$  Therefore, the dimensions given for all the orbits are correct.

As in the above cases, it is very easy to show that each orbit is contained in the locally closed set $S$ in the second column by simple matrix calculations, so we leave this part of the proof to the reader as we did in type $A_{3} .$  As in the previous case, we will provide, for each of the $61$ orbits, an element $b$ in $\mathfrak{b}$ which has the property that $b .n$ is an arbitrary element of this closed set $S$, thereby showing the reverse containments.  This will take some space as for each orbit, treated in a separate table, we give the orbit representative in the first row, the locally closed set $S$ in the second row, the form $m$ of an arbitrary element of $S$ as an ordered $10$-tuple in $\mathfrak{n} \approx k^{10}$ in the third row, and the last row (or occasionally last two rows...) will contain the desired element of $b ,$ as usual, written as a product of an element of $T$ with several elements from root groups.  In each case, the reader should actually check that \thinspace $b .n =m .$  This will prove that all the defining equations are correct for each orbit.

\begin{center}
\begin{tabular}[c]{|c|}\hline
$0$ \\
\hline
$Z(X_{1} ,X_{2} ,X_{3} ,X_{4} ,X_{12} ,X_{23} ,X_{34} ,X_{13} ,X_{24} ,X_{14})$ \\
\hline
$m =(0 ,0 ,0 ,0 ,0 ,0 ,0 ,0 ,0 ,0)$ \\
\hline
any $b \in B$ \\
\hline
\end{tabular}\end{center}\par
\begin{center}
\begin{tabular}[c]{|c|}\hline
$x_{14}$ \\
\hline
$Z(X_{1} ,X_{2} ,X_{3} ,X_{4} ,X_{12} ,X_{23} ,X_{34} ,X_{13} ,X_{24}) \cap V(X_{14})$ \\
\hline
$m =(0 ,0 ,0 ,0 ,0 ,0 ,0 ,0 ,0 ,z) ,\ z \neq 0$ \\
\hline
$b =T(1 ,1 ,1 ,z)$ \\
\hline
\end{tabular}\end{center}\par
\begin{center}
\begin{tabular}[c]{|c|}\hline
$x_{13}$ \\
\hline
$Z(X_{1} ,X_{2} ,X_{3} ,X_{4} ,X_{12} ,X_{23} ,X_{34} ,X_{24}) \cap V(X_{13})$ \\
\hline
$m =(0 ,0 ,0 ,0 ,0 ,0 ,0 ,x ,0 ,z) ,\ x \neq 0$ \\
\hline
$b =T(x ,1 ,1 ,1)U_{4}\genfrac{(}{)}{}{}{z}{x^{2}}$ \\
\hline
\end{tabular}\end{center}\par
\begin{center}
\begin{tabular}[c]{|c|}\hline
$x_{24}$ \\
\hline
$Z(X_{1} ,X_{2} ,X_{3} ,X_{4} ,X_{12} ,X_{23} ,X_{34} ,X_{13}) \cap V(X_{24})$ \\
\hline
$m =(0 ,0 ,0 ,0 ,0 ,0 ,0 ,0 ,y ,z) ,\ y \neq 0$ \\
\hline
$b =T(y ,1 ,1 ,z)U_{1}\genfrac{(}{)}{}{}{z}{y^{2}}$ \\
\hline
\end{tabular}\end{center}\par
\begin{center}
\begin{tabular}[c]{|c|}\hline
$x_{13} +x_{24}$ \\
\hline
$Z(X_{1} ,X_{2} ,X_{3} ,X_{4} ,X_{12} ,X_{23} ,X_{34}) \cap V(X_{13} ,X_{24})$ \\
\hline
$m =(0 ,0 ,0 ,0 ,0 ,0 ,0 ,x ,y ,z) ,\ x ,y \neq 0$ \\
\hline
$b =T(x ,1 ,\frac{y}{x} ,1)U_{1}\genfrac{(}{)}{}{}{z}{xy}$ \\
\hline
\end{tabular}\end{center}\par
\begin{center}
\begin{tabular}[c]{|c|}\hline
$x_{12}$ \\
\hline
$Z(X_{1} ,X_{2} ,X_{3} ,X_{4} ,X_{23} ,X_{34} ,X_{24}) \cap V(X_{12})$ \\
\hline
$m =(0 ,0 ,0 ,0 ,u ,0 ,0 ,x ,0 ,z) ,\ u \neq 0$ \\
\hline
$b =T(u ,1 ,1 , -1)U_{3}\genfrac{(}{)}{}{}{x}{u}U_{34}\genfrac{(}{)}{}{}{z}{u^{2}}$ \\
\hline
\end{tabular}\end{center}\par
\begin{center}
\begin{tabular}[c]{|c|}\hline
$x_{23}$ \\
\hline
$Z(X_{1} ,X_{2} ,X_{3} ,X_{4} ,X_{12} ,X_{34} ,X_{13}X_{24} -X_{23}X_{14}) \cap V(X_{23})$ \\
\hline
$m =(0 ,0 ,0 ,0 ,0 ,v ,0 ,x ,y ,z) ,\ v \neq 0 ,\ xy -vz =0$ \\
\hline
$b =T(1 ,v , -1 ,1)U_{1}\left (x\right )U_{4}\genfrac{(}{)}{}{}{y}{v^{2}}$ \\
\hline
\end{tabular}\end{center}\par
\begin{center}
\begin{tabular}[c]{|c|}\hline
$x_{34}$ \\
\hline
$Z(X_{1} ,X_{2} ,X_{3} ,X_{4} ,X_{12} ,X_{23} ,X_{13}) \cap V(X_{34})$ \\
\hline
$m =(0 ,0 ,0 ,0 ,0 ,0 ,w ,0 ,y ,z) ,\ w \neq 0$ \\
\hline
$b =T(w ,1 ,1 ,1)U_{2}\genfrac{(}{)}{}{}{y}{w}U_{23}\genfrac{(}{)}{}{}{z}{w^{2}}$ \\
\hline
\end{tabular}\end{center}\par
\begin{center}
\begin{tabular}[c]{|c|}\hline
$x_{12} +x_{24}$ \\
\hline
$Z(X_{1} ,X_{2} ,X_{3} ,X_{4} ,X_{23} ,X_{34}) \cap V(X_{12} ,X_{24})$ \\
\hline
$m =(0 ,0 ,0 ,0 ,u ,0 ,0 ,x ,y ,z) ,\ u ,y \neq 0$ \\
\hline
$b =T(u ,\sqrt{\frac{y}{u}} ,1 ,1)U_{1}\genfrac{(}{)}{}{}{z}{\sqrt{u^{3}y}}U_{3}\genfrac{(}{)}{}{}{ -x}{u}$ \\
\hline
\end{tabular}\end{center}\par
\begin{center}
\begin{tabular}[c]{|c|}\hline
$x_{23} +x_{14}$ \\
\hline
$Z(X_{1} ,X_{2} ,X_{3} ,X_{4} ,X_{12} ,X_{34}) \cap V(X_{23 ,\thinspace }X_{13}X_{24} -X_{23}X_{14})$ \\
\hline
$m =(0 ,0 ,0 ,0 ,0 ,v ,0 ,x ,y ,z) ,\ v ,xy -vz \neq 0$ \\
\hline
$b =T(1 ,v ,\frac{vz -xy}{v^{2}} ,1)U_{1}\left .(x\right )U_{4}\genfrac{(}{)}{}{}{y}{xy -vz}$ \\
\hline
\end{tabular}\end{center}\par
\begin{center}
\begin{tabular}[c]{|c|}\hline
$x_{34} +x_{13}$ \\
\hline
$Z(X_{1} ,X_{2} ,X_{3} ,X_{4} ,X_{12} ,X_{23}) \cap V(X_{34} ,X_{13})$ \\
\hline
$m =(0 ,0 ,0 ,0 ,0 ,0 ,w ,x ,y ,z) ,\ w ,x \neq 0$ \\
\hline
$b =T(\sqrt{wx} ,1 ,1 ,\sqrt{\frac{w}{x}})U_{4}\genfrac{(}{)}{}{}{ -z}{\sqrt{w^{3}x}}U_{2}\genfrac{(}{)}{}{}{y}{w}$ \\
\hline
\end{tabular}\end{center}\par
\begin{center}
\begin{tabular}[c]{|c|}\hline
$x_{1}$ \\
\hline
$Z(X_{2} ,X_{3} ,X_{4} ,X_{23} ,X_{34} ,X_{24}) \cap V(X_{1})$ \\
\hline
$m =(q ,0 ,0 ,0 ,u ,0 ,0 ,x ,0 ,z) ,\ q \neq 0$ \\
\hline
$T(q ,1 ,1 ,1)U_{2}\genfrac{(}{)}{}{}{ -u}{q}U_{23}\genfrac{(}{)}{}{}{ -x}{q}U_{24}\genfrac{(}{)}{}{}{ -z}{q^{2}}$ \\
\hline
\end{tabular}\end{center}\par
\begin{center}
\begin{tabular}[c]{|c|}\hline
$x_{2}$ \\
\hline
$Z(X_{1} ,X_{3} ,X_{4} ,X_{34} ,\thinspace X_{2}X_{13} -X_{12}X_{23} ,\thinspace X_{2}X_{14} -X_{12}X_{24}) \cap V(X_{2})$ \\
\hline
$m =(0 ,r ,0 ,0 ,u ,v ,0 ,x ,y ,z) ,\ \thinspace rx -uv ,rz -yu =0 ,\thinspace r \neq 0$ \\
\hline
$T(1 ,r ,1 ,1)U_{1}(u)U_{3}\genfrac{(}{)}{}{}{ -v}{r}U_{34}\genfrac{(}{)}{}{}{ -y}{r^{2}}$ \\
\hline
\end{tabular}\end{center}\par
\begin{center}
\begin{tabular}[c]{|c|}\hline
$x_{3}$ \\
\hline
$Z(X_{1} ,X_{2} ,X_{4} ,X_{12} ,\thinspace X_{3}X_{24} -X_{23}X_{34} ,\thinspace X_{3}X_{14} -X_{13}X_{34}) \cap V(X_{3})$ \\
\hline
$m =(0 ,0 ,s ,0 ,0 ,v ,w ,x ,y ,z) ,\ \ sy -vw ,sz -wx =0 ,\thinspace s \neq 0$ \\
\hline
$b =T(s ,1 ,1 ,\frac{1}{s})U_{4}\left .( -w\right )U_{2}\genfrac{(}{)}{}{}{v}{s}U_{12}\genfrac{(}{)}{}{}{x}{s^{2}}$ \\
\hline
\end{tabular}\end{center}\par
\begin{center}
\begin{tabular}[c]{|c|}\hline
$x_{4}$ \\
\hline
$Z(X_{1} ,X_{2} ,X_{3} ,X_{12} ,X_{23} ,X_{13}) \cap V(X_{4})$ \\
\hline
$m =(0 ,0 ,0 ,t ,0 ,0 ,w ,0 ,y ,z) ,\ t \neq 0$ \\
\hline
$T(t ,1 ,1 ,1)U_{3}\genfrac{(}{)}{}{}{w}{t}U_{23}\genfrac{(}{)}{}{}{y}{t}U_{13}\genfrac{(}{)}{}{}{z}{t^{2}}$ \\
\hline
\end{tabular}\end{center}\par
\begin{center}
\begin{tabular}[c]{|c|}\hline
$x_{12} +x_{23}$ \\
\hline
$Z(X_{1} ,X_{2} ,X_{3} ,X_{4} ,X_{34}) \cap V(X_{12} ,X_{23})$ \\
\hline
$m =(0 ,0 ,0 ,0 ,u ,v ,0 ,x ,y ,z) ,\ u ,v \neq 0$ \\
\hline
$T(u ,v ,1 ,1)U_{1}\genfrac{(}{)}{}{}{x}{u}U_{4}\genfrac{(}{)}{}{}{ -y}{uv^{2}}U_{34}\genfrac{(}{)}{}{}{xy -zv}{u^{2}v^{2}}$ \\
\hline
\end{tabular}\end{center}\par
\begin{center}
\begin{tabular}[c]{|c|}\hline
$x_{12} +x_{34}$ \\
\hline
$Z(X_{1} ,X_{2} ,X_{3} ,X_{4} ,X_{23}) \cap V(X_{12} ,X_{34})$ \\
\hline
$m =(0 ,0 ,0 ,0 ,u ,0 ,w ,x ,y ,z) ,\ u ,w \neq 0$ \\
\hline
$T(u ,\frac{u}{w} ,1 ,1)U_{3}\genfrac{(}{)}{}{}{ -x}{u}U_{2}\genfrac{(}{)}{}{}{yu}{w^{2}}U_{12}\genfrac{(}{)}{}{}{z}{uw}$ \\
\hline
\end{tabular}\end{center}\par
\begin{center}
\begin{tabular}[c]{|c|}\hline
$x_{23} +x_{34}$ \\
\hline
$Z(X_{1} ,X_{2} ,X_{3} ,X_{4} ,X_{12}) \cap V(X_{23} ,X_{34})$ \\
\hline
$m =(0 ,0 ,0 ,0 ,u ,0 ,w ,x ,y ,z) ,\ v ,w \neq 0$ \\
\hline
$T(vw ,1 ,1 ,\frac{1}{v})U_{4}\genfrac{(}{)}{}{}{ -y}{w}U_{1}\genfrac{(}{)}{}{}{x}{v^{2}w}U_{12}\genfrac{(}{)}{}{}{zv -xy}{v^{2}w^{2}}$ \\
\hline
\end{tabular}\end{center}\par
\begin{center}
\begin{tabular}[c]{|c|}\hline
$x_{1}$ \\
\hline
$Z(X_{2} ,X_{3} ,X_{4} ,X_{23} ,X_{34} ,X_{24}) \cap V(X_{1})$ \\
\hline
$m =(q ,0 ,0 ,0 ,u ,0 ,0 ,x ,0 ,z) ,\ q \neq 0$ \\
\hline
$T(q ,1 ,1 ,1)U_{2}\genfrac{(}{)}{}{}{ -u}{q}U_{23}\genfrac{(}{)}{}{}{ -x}{q}U_{24}\genfrac{(}{)}{}{}{ -z}{q^{2}}$ \\
\hline
\end{tabular}\end{center}\par
\begin{center}
\begin{tabular}[c]{|c|}\hline
$x_{1} +x_{24}$ \\
\hline
$Z(X_{2} ,X_{3} ,X_{4} ,X_{23} ,X_{34}) \cap V(X_{1} ,X_{24})$ \\
\hline
$m =(q ,0 ,0 ,0 ,u ,0 ,0 ,x ,y ,z) ,\ q ,y \neq 0$ \\
\hline
$T(q^{\frac{2}{3}}y^{\frac{1}{3}} ,y^{\frac{1}{3}}q^{ -\frac{1}{3}} ,1 ,1)U_{2}\genfrac{(}{)}{}{}{ -u}{q^{\frac{2}{3}}y^{\frac{1}{3}}}U_{23}\genfrac{(}{)}{}{}{ -x}{q^{\frac{2}{3}}y^{\frac{1}{3}}}U_{24}\genfrac{(}{)}{}{}{ -z}{qy}$ \\
\hline
\end{tabular}\end{center}\par
\begin{center}
\begin{tabular}[c]{|c|}\hline
$x_{1} +x_{34}$ \\
\hline
$Z(X_{2} ,X_{3} ,X_{4} ,X_{23} ,X_{12}X_{34} +X_{1}X_{24}) \cap V(X_{1} ,X_{34})$ \\
\hline
$m =(q ,0 ,0 ,0 ,u ,0 ,w ,x ,y ,z) ,\ uw +qy =0 ,\thinspace q ,w \neq 0$ \\
\hline
$T(\sqrt{qw} ,\sqrt{\frac{w}{q}} ,1 ,1)U_{2}\genfrac{(}{)}{}{}{ -u}{\sqrt{qw}}U_{23}\genfrac{(}{)}{}{}{ -x}{\sqrt{qw}}U_{24}\genfrac{(}{)}{}{}{ -z}{\sqrt{qw^{3}}}$ \\
\hline
\end{tabular}\end{center}\par
\begin{center}
\begin{tabular}[c]{|c|}\hline
$x_{2} +x_{14}$ \\
\hline
$Z(X_{1} ,X_{3} ,X_{4} ,X_{34} ,\thinspace X_{12}X_{23} -X_{2}X_{13}) \cap V(X_{2} ,X_{12}X_{24} -X_{2}X_{14})$ \\
\hline
$m =(0 ,r ,0 ,0 ,u ,v ,0 ,x ,y ,z) ,\ uv -rx =0 ,\thinspace r ,yu -rz \neq 0$ \\
\hline
$T(\frac{\sqrt{rz -yu}}{r} ,r ,1 ,1)U_{1}\genfrac{(}{)}{}{}{ur}{\sqrt{rz -yu}}U_{3}\genfrac{(}{)}{}{}{ -v}{r}U_{34}\genfrac{(}{)}{}{}{ -y}{r\sqrt{rz -yu}}$ \\
\hline
\end{tabular}\end{center}\par
\begin{center}
\begin{tabular}[c]{|c|}\hline
$x_{3} +x_{14}$ \\
\hline
$Z(X_{1} ,X_{2} ,X_{4} ,X_{12} ,\thinspace X_{23}X_{34} -X_{3}X_{24}) \cap V(X_{3} ,\thinspace X_{34}X_{13} -X_{3}X_{14})$ \\
\hline
$m =(0 ,0 ,s ,0 ,0 ,v ,w ,x ,y ,z) ,\ vw -sy =0 ,\ s ,wx -sz \neq 0$ \\
\hline
$T(\sqrt{sz -wx} ,1 ,1 ,\frac{1}{s})U_{4}\genfrac{(}{)}{}{}{ -ws}{\sqrt{sz -wx}}U_{2}\genfrac{(}{)}{}{}{v}{s}U_{12}\genfrac{(}{)}{}{}{x}{s\sqrt{sz -wx}}$ \\
\hline
\end{tabular}\end{center}\par
\begin{center}
\begin{tabular}[c]{|c|}\hline
$x_{4} +x_{12}$ \\
\hline
$Z(X_{1} ,X_{2} ,X_{3} ,X_{23} ,X_{12}X_{34} +X_{4}X_{13}) \cap V(X_{4} ,X_{12})$ \\
\hline
$m =(0 ,0 ,0 ,t ,u ,0 ,w ,x ,y ,z) ,\ uw +tx =0 ,\ t ,u \neq 0$ \\
\hline
$T(u ,1 ,1 ,\sqrt{\frac{t}{u}})U_{3}\genfrac{(}{)}{}{}{w}{\sqrt{tu}}U_{23}\genfrac{(}{)}{}{}{y}{\sqrt{tu}}U_{13}\genfrac{(}{)}{}{}{z}{\sqrt{tu^{3}}}$ \\
\hline
\end{tabular}\end{center}\par
\begin{center}
\begin{tabular}[c]{|c|}\hline
$x_{4} +x_{13}$ \\
\hline
$Z(X_{1} ,X_{2} ,X_{3} ,X_{12} ,X_{23}) \cap V(X_{4} ,X_{13})$ \\
\hline
$m =(0 ,0 ,0 ,t ,0 ,0 ,w ,x ,y ,z) ,\ t ,x \neq 0$ \\
\hline
$T(t^{\frac{1}{3}}x^{\frac{2}{3}} ,1 ,1 ,t^{\frac{1}{3}}x^{ -\frac{1}{3}})U_{3}\genfrac{(}{)}{}{}{w}{t^{\frac{2}{3}}x^{\frac{1}{3}}}U_{23}\genfrac{(}{)}{}{}{y}{t^{\frac{2}{3}}x^{\frac{1}{3}}}U_{13}\genfrac{(}{)}{}{}{z}{tx}$ \\
\hline
\end{tabular}\end{center}\par
\begin{center}
\begin{tabular}[c]{|c|}\hline
$x_{12} +x_{23} +x_{34}$ \\
\hline
$Z(X_{1} ,X_{2} ,X_{3} ,X_{4}) \cap V(X_{2} ,X_{23} ,X_{34})$ \\
\hline
$m =(0 ,0 ,0 ,0 ,u ,v ,w ,x ,y ,z) ,\ u ,v ,w \neq 0$ \\
\hline
$T(u ,\sqrt{\frac{vw}{u}} ,1 ,\sqrt{\frac{w}{uv}})U_{1}\genfrac{(}{)}{}{}{x\sqrt{w}}{\sqrt{u^{3}v}}U_{2}\genfrac{(}{)}{}{}{y\sqrt{u}}{\sqrt{vw^{3}}}U_{12}\genfrac{(}{)}{}{}{vz -xy}{uvw}$ \\
\hline
\end{tabular}\end{center}\par
\begin{center}
\begin{tabular}[c]{|c|}\hline
$x_{1} +x_{23}$ \\
\hline
$Z(X_{2} ,X_{3} ,X_{4} ,X_{34}) \cap V(X_{1} ,X_{23})$ \\
\hline
$m =(q ,0 ,0 ,0 ,u ,v ,0 ,x ,y ,z) ,\ q ,v \neq 0$ \\
\hline
$T(qv ,v ,1 ,1)U_{1}\genfrac{(}{)}{}{}{x}{qv}U_{2}\genfrac{(}{)}{}{}{ -u}{qv}U_{4}(\genfrac{(}{)}{}{}{ -y}{qv^{3}})U_{24}.\genfrac{(}{)}{}{}{xy -vz}{q^{2}v^{4}}$ \\
\hline
\end{tabular}\end{center}\par
\begin{center}
\begin{tabular}[c]{|c|}\hline
$x_{1} +x_{34} +x_{24}$ \\
\hline
$Z(X_{2} ,X_{3} ,X_{4} ,X_{23}) \cap V(X_{1} ,X_{34} ,\thinspace X_{1}X_{24} +X_{12}X_{34})$ \\
\hline
$m =(q ,0 ,0 ,0 ,u ,0 ,w ,x ,y ,z) ,\ q ,w ,qy +uw \neq 0$ \\
\hline
$T(q ,1 ,\frac{qw}{qy +uw} ,\frac{(qy +uw)^{2}}{q^{3}w})U_{2}\genfrac{(}{)}{}{}{ -uw}{qy +uw}U_{23}\genfrac{(}{)}{}{}{ -x(qy +uw)^{2}}{q^{4}w}U_{24}\genfrac{(}{)}{}{}{ -z}{qy +uw}$ \\
\hline
\end{tabular}\end{center}\par
\begin{center}
\begin{tabular}[c]{|c|}\hline
$x_{2} +x_{13}$ \\
\hline
$Z(X_{1} ,X_{3} ,X_{4} ,X_{23} ,X_{34}) \cap V(X_{2} ,\thinspace X_{12}X_{23} -X_{2}X_{13})$ \\
\hline
$m =(0 ,r ,0 ,0 ,u ,v ,0 ,x ,y ,z) ,\ r ,uv -rx \neq 0$ \\
\hline
$T(\frac{rx -uv}{r} ,r ,1 ,1)U_{1}\genfrac{(}{)}{}{}{ur}{rx -uv}U_{3}\genfrac{(}{)}{}{}{ -v}{r}U_{34}\genfrac{(}{)}{}{}{y}{r(uv -rx)}U_{4}\genfrac{(}{)}{}{}{yu -rz}{(rx -uv)^{2}}$ \\
\hline
\end{tabular}\end{center}\par
\begin{center}
\begin{tabular}[c]{|c|}\hline
$x_{2} +x_{34}$ \\
\hline
$Z(X_{1} ,X_{3} ,X_{4} ,X_{2}X_{13} -X_{12}X_{23}) \cap V(X_{2} ,X_{34})$ \\
\hline
$m =(0 ,r ,0 ,0 ,u ,v ,w ,x ,y ,z) ,\ rx -uv =0 ,\ r ,w \neq 0$ \\
\hline
$T(\frac{w}{r} ,r ,1 ,1)U_{3}\genfrac{(}{)}{}{}{ -v}{r}U_{1}\genfrac{(}{)}{}{}{ru}{w}U_{2}\genfrac{(}{)}{}{}{y}{rw}U_{12}\genfrac{(}{)}{}{}{rz -yu}{w^{2}}$ \\
\hline
\end{tabular}\end{center}\par
\begin{center}
\begin{tabular}[c]{|c|}\hline
$x_{3} +x_{12}$ \\
\hline
$Z(X_{1} ,X_{2} ,X_{4} ,X_{3}X_{24} -X_{23}X_{34}) \cap V(X_{3} ,X_{12})$ \\
\hline
$m =(0 ,0 ,s ,0 ,u ,v ,w ,x ,y ,z) ,\ sy -vw =0 ,\ s ,u \neq 0$ \\
\hline
$T(u ,1 ,1 ,\frac{1}{s})U_{2}\genfrac{(}{)}{}{}{v}{s}U_{4}\genfrac{(}{)}{}{}{ -sw}{u}U_{3}\genfrac{(}{)}{}{}{ -x}{su}U_{34}\genfrac{(}{)}{}{}{wx -sz}{u^{2}}$ \\
\hline
\end{tabular}\end{center}\par
\begin{center}
\begin{tabular}[c]{|c|}\hline
$x_{3} +x_{24}$ \\
\hline
$Z(X_{1} ,X_{2} ,X_{4} ,X_{12}) \cap V(X_{3} ,\thinspace X_{23}X_{34} -X_{3}X_{24})$ \\
\hline
$m =(0 ,0 ,s ,0 ,0 ,v ,w ,x ,y ,z) ,\ s ,vw -sy \neq 0$ \\
\hline
$T(sy -vw ,1 ,1 ,\frac{1}{s})U_{4}\genfrac{(}{)}{}{}{ws}{vw -sy}U_{2}\genfrac{(}{)}{}{}{v}{s}U_{12}\genfrac{(}{)}{}{}{x}{s(sy -vw)}U_{1}\genfrac{(}{)}{}{}{sz -xw}{(sy -vw)^{2}}$ \\
\hline
\end{tabular}\end{center}\par
\begin{center}
\begin{tabular}[c]{|c|}\hline
$x_{4} +x_{12} +x_{13}$ \\
\hline
$Z(X_{1} ,X_{2} ,X_{3} ,X_{23}) \cap V(X_{4} ,X_{12} ,\thinspace X_{4}X_{13} +X_{12}X_{34})$ \\
\hline
$m =(0 ,0 ,0 ,t ,u ,0 ,w ,x ,y ,z) ,\ t ,u ,tx +uw \neq 0$ \\
\hline
$T(u ,\frac{(tx +uw)^{2}}{tu^{3}} ,1 ,\frac{tu}{tx +uw})U_{3}\genfrac{(}{)}{}{}{uw}{tx +uw}U_{23}\genfrac{(}{)}{}{}{ytu^{4}}{(tx +uw)^{3}}U_{13}\genfrac{(}{)}{}{}{z}{tx +uw}$ \\
\hline
\end{tabular}\end{center}\par
\begin{center}
\begin{tabular}[c]{|c|}\hline
$x_{1} +x_{3}$ \\
\hline
$Z(X_{2} ,X_{4} ,\thinspace X_{1}X_{23} +X_{3}X_{12} ,\thinspace X_{3}X_{24} -X_{23}X_{34}) \cap V(X_{1} ,\thinspace X_{3})$ \\
\hline
$m =(q ,0 ,s ,0 ,u ,v ,w ,x ,y ,z) ,\ qv +su ,sy -vw =0 ,\ q ,s \neq 0$ \\
\hline
$T(q ,1 ,s ,1)U_{4}\genfrac{(}{)}{}{}{ -w}{qs^{2}}U_{2}\left (v\right )U_{23}\genfrac{(}{)}{}{}{ -x}{q}U_{24}\genfrac{(}{)}{}{}{wx -sz}{(qs)^{2}}$ \\
\hline
\end{tabular}\end{center}\par
\begin{center}
\begin{tabular}[c]{|c|}\hline
$x_{1} +x_{4}$ \\
\hline
$Z(X_{2} ,X_{3} ,X_{23} ,X_{1}X_{24} +X_{12}X_{34} +X_{4}X_{13}) \cap V(X_{1} ,X_{4})$ \\
\hline
$m =(q ,0 ,0 ,t ,u ,0 ,w ,x ,y ,z) ,\ qy +uw +tx =0 ,\ q ,t \neq 0$ \\
\hline
$T(q ,1 ,\frac{t}{q} ,1)U_{2}\genfrac{(}{)}{}{}{ -tu}{q^{2}}U_{3}\genfrac{(}{)}{}{}{qw}{t^{2}}U_{23}\genfrac{(}{)}{}{}{ -x}{q}U_{24}\genfrac{(}{)}{}{}{ -z}{qt}$ \\
\hline
\end{tabular}\end{center}\par
\begin{center}
\begin{tabular}[c]{|c|}\hline
$x_{2} +x_{4}$ \\
\hline
$Z(X_{1} ,X_{3} ,X_{4}X_{23} +X_{2}X_{34} ,X_{2}X_{13} -X_{12}X_{23}) \cap V(X_{2} ,X_{4})$ \\
\hline
$m =(0 ,r ,0 ,t ,u ,v ,w ,x ,y ,z) ,\ tv +rw ,rx -uv =0 ,\ r ,t \neq 0$ \\
\hline
$T(rt ,1 ,\frac{1}{r} ,1)U_{1}\genfrac{(}{)}{}{}{u}{tr^{2}}U_{3}\left ( -v\right )U_{23}genfrac{(}{)}{}{}{y}{t}U_{13}\genfrac{(}{)}{}{}{rz -uy}{(tr)^{2}}$ \\
\hline
\end{tabular}\end{center}\par
\begin{center}
\begin{tabular}[c]{|c|}\hline
$x_{1} +x_{23} +x_{34}$ \\
\hline
$Z(X_{2} ,X_{3} ,X_{4}) \cap V(X_{1} ,X_{23} ,X_{34})$ \\
\hline
$m =(q ,0 ,0 ,0 ,u ,v ,w ,x ,y ,z) ,\ q ,v ,w \neq 0$ \\
\hline
{\begin{tabular}[c]{c}$T(q ,1 ,\sqrt{\frac{vw}{q}} ,\frac{1}{v})$$ \cdot $ \\
$U_{4}\genfrac{(}{)}{}{}{ -\sqrt{v}(qy +uw)}{\sqrt{q^{3}w}}U_{2}\genfrac{(}{)}{}{}{ -u\sqrt{vw}}{\sqrt{q^{3}}}U_{23}genfrac{(}{)}{}{}{ -x}{qv}).U_{24}\genfrac{(}{)}{}{}{qxy -qvz +uwx}{\sqrt{q^{5}vw}}$
\end{tabular}} \\
\hline
\end{tabular}\end{center}\par
\begin{center}
\begin{tabular}[c]{|c|}\hline
$x_{2} +x_{34} +x_{13}$ \\
\hline
$Z(X_{1} ,X_{3} ,X_{4}) \cap V(X_{2} ,X_{34} ,\thinspace X_{2}X_{13} -X_{12}X_{23})$ \\
\hline
$m =(0 ,r ,0 ,0 ,u ,v ,w ,x ,y ,z) ,\ r ,w ,rx -uv \neq 0$ \\
\hline
{\begin{tabular}[c]{c}$T(\frac{\sqrt{w(rx -uv)}}{r} ,r ,1 ,\sqrt{\frac{w}{rx -uv}})$$ \cdot $ \\
$U_{1}\genfrac{(}{)}{}{}{ru}{\sqrt{w(rx -uv)}}U_{2}\genfrac{(}{)}{}{}{y}{rw}U_{12}\genfrac{(}{)}{}{}{rz -yu}{\sqrt{w^{3}(rx -uv)}}U_{3}\genfrac{(}{)}{}{}{ -v\sqrt{w}}{r\sqrt{rx -uv}}$
\end{tabular}} \\
\hline
\end{tabular}\end{center}\par
\begin{center}
\begin{tabular}[c]{|c|}\hline
$x_{3} +x_{12} +x_{24}$ \\
\hline
$Z(X_{1} ,X_{2} ,X_{4}) \cap V(X_{3} ,X_{12} ,X_{3}X_{24} -X_{23}X_{34})$ \\
\hline
$m =(0 ,0 ,s ,0 ,u ,v ,w ,x ,y ,z) ,\ s ,u ,sy -vw \neq 0$ \\
\hline
{\begin{tabular}[c]{c}$T(u ,\sqrt{\frac{sy -vw}{u}} ,1 ,\frac{1}{s}) \cdot $ \\
$U_{4}\genfrac{(}{)}{}{}{ -ws}{\sqrt{u(sy -vw)}}U_{3}\genfrac{(}{)}{}{}{ -x}{su}U_{34}\genfrac{(}{)}{}{}{wx -sz}{\sqrt{u^{3}(sy -vw)}}U_{2}\genfrac{(}{)}{}{}{v\sqrt{u}}{s\sqrt{sy -vw}}$
\end{tabular}} \\
\hline
\end{tabular}\end{center}\par
\begin{center}
\begin{tabular}[c]{|c|}\hline
$x_{4} +x_{12} +x_{23}$ \\
\hline
$Z(X_{1} ,X_{2} ,X_{3}) \cap V(X_{4} ,X_{12} ,X_{23})$ \\
\hline
$m =(0 ,0 ,0 ,t ,u ,v ,w ,x ,y ,z) ,\ t ,v ,u \neq 0$ \\
\hline
$T(\sqrt{\frac{tu}{v}} ,v ,\sqrt{\frac{t}{uv}} ,1)U_{1}\genfrac{(}{)}{}{}{(tx +uw)\sqrt{v}}{\sqrt{t^{3}u}}U_{3}\genfrac{(}{)}{}{}{w\sqrt{uv}}{\sqrt{t^{3}}}U_{23}\genfrac{(}{)}{}{}{y}{tv}U_{13}\genfrac{(}{)}{}{}{tvz -txy -uwy}{\sqrt{t^{5}uv}}$ \\
\hline
\end{tabular}\end{center}\par
\begin{center}
\begin{tabular}[c]{|c|}\hline
$x_{1} +x_{2}$ \\
\hline
$Z(X_{3} ,X_{4} ,X_{34}) \cap V(X_{1} ,X_{2})$ \\
\hline
$m =(q ,r ,0 ,0 ,u ,v ,0 ,x ,y ,z) ,\ q ,r \neq 0$ \\
\hline
{\begin{tabular}[c]{c}$T(q^{\frac{2}{3}}r^{\frac{1}{3}} ,r^{\frac{1}{3}}q^{ -\frac{1}{3}} ,r^{ -\frac{2}{3}}q^{ -\frac{1}{3}} ,1)$ $ \cdot $ \\
$U_{1}\genfrac{(}{)}{}{}{u}{qr}U_{3}\genfrac{(}{)}{}{}{ -vq^{\frac{1}{3}}}{r^{\frac{1}{3}}}U_{23}\genfrac{(}{)}{}{}{uv -rx}{q^{\frac{2}{3}}r^{\frac{4}{3}}}U_{34}\genfrac{(}{)}{}{}{ -yq^{\frac{1}{3}}}{r^{\frac{1}{3}}}U_{24}\genfrac{(}{)}{}{}{yu -rz}{q^{\frac{2}{3}}r^{\frac{4}{3}}}$
\end{tabular}} \\
\hline
\end{tabular}\end{center}\par
\begin{center}
\begin{tabular}[c]{|c|}\hline
$x_{1} +x_{3} +x_{12}$ \\
\hline
$Z(X_{2} ,X_{4} ,\thinspace X_{3}X_{23} -X_{23}X_{34}) \cap V(X_{1} ,X_{3}\thinspace X_{1}X_{23} +X_{3}X_{12})$ \\
\hline
$m =(q ,0 ,s ,0 ,u ,v ,w ,x ,y ,z) ,\ sy -vw =0 ,\ q ,s ,qv +us \neq 0$ \\
\hline
{\begin{tabular}[c]{c}$T(q ,1 ,\frac{qs}{qv +us} ,\frac{q}{qv +us})$ \\
$U_{2}\genfrac{(}{)}{}{}{qv}{qv +us}U_{4}\genfrac{(}{)}{}{}{ -w(qv +us)^{3}}{q^{4}s^{2}}U_{12}\genfrac{(}{)}{}{}{x}{qv +us}U_{24}\genfrac{(}{)}{}{}{(qv +us)^{2}(wx -sz)}{q^{4}s^{2}}$
\end{tabular}} \\
\hline
\end{tabular}\end{center}\par
\begin{center}
\begin{tabular}[c]{|c|}\hline
$x_{1} +x_{3} +x_{24}$ \\
\hline
$Z(X_{2} ,X_{4} ,\thinspace X_{1}X_{23} +X_{3}X_{12}) \cap V(X_{1} ,X_{3} ,X_{3}X_{24} -X_{23}X_{34})$ \\
\hline
$m =(q ,0 ,s ,0 ,u ,v ,w ,x ,y ,z) ,\ qv +su =0 ,\ q ,s .sy -vw \neq 0$ \\
\hline
{\begin{tabular}[c]{c}$T(\frac{q\frac{2}{3}(sy -vw)^{\frac{1}{3}}}{s^{\frac{2}{3}}} ,\frac{(sy -vw)^{\frac{1}{3}}}{q^{\frac{1}{3}}s^{\frac{2}{3}}} ,s ,1) \cdot $ \\
$U_{1}\genfrac{(}{)}{}{}{zs -wx}{q(sy -vw)}U_{2}\genfrac{(}{)}{}{}{vq^{\frac{1}{3}}s^{\frac{2}{3}}}{(sy -vw)^{\frac{1}{3}}}U_{4}\genfrac{(}{)}{}{}{ -w}{q^{\frac{1}{3}}s^{\frac{2}{3}}(sy -vw)^{\frac{2}{3}}}U_{23}\genfrac{(}{)}{}{}{(vz -xy)s^{\frac{5}{3}}}{q^{\frac{2}{3}}(sy -vw)^{\frac{4}{3}}}$
\end{tabular}} \\
\hline
\end{tabular}\end{center}\par
\begin{center}
\begin{tabular}[c]{|c|}\hline
$x_{2} +x_{3}$ \\
\hline
$Z(X_{1} ,X_{4} ,\thinspace X_{3}X_{12}X_{24} +X_{2}X_{34}X_{13} -X_{12}X_{23}X_{34} -X_{2}X_{3}X_{14}) \cap V(X_{2} ,X_{3})$ \\
\hline
$m =(0 ,r ,s ,0 ,u ,v ,w ,x ,y ,z) ,\ syu +rwx -uvw -rsz =0 ,\ r ,s \neq 0$ \\
\hline
$T(\frac{s}{r} ,r ,1 ,\frac{1}{s})U_{1}\genfrac{(}{)}{}{}{ur}{s}U_{2}\genfrac{(}{)}{}{}{v}{rs}U_{12}\genfrac{(}{)}{}{}{rx -uv}{s^{2}}U_{4}( -w)U_{34}\genfrac{(}{)}{}{}{vw -sy}{rs}$ \\
\hline
\end{tabular}\end{center}\par
\begin{center}
\begin{tabular}[c]{|c|}\hline
$x_{1} +x_{4} +x_{24}$ \\
\hline
$Z(X_{2} ,X_{3} ,X_{23}) \cap V(X_{1} ,X_{4} ,\thinspace X_{1}X_{24} +X_{12}X_{44} +X_{4}X_{13})$ \\
\hline
$m =(q ,0 ,0 ,t ,u ,0 ,w ,x ,y ,z) ,\ q ,t ,qy +uw +tx \neq 0$ \\
\hline
{\begin{tabular}[c]{c}$T(q ,1 ,\frac{(qy +uw +tx)^{2}}{q^{3}t} ,\frac{qt}{qy +uw +tx}) \cdot $ \\
$U_{2}\genfrac{(}{)}{}{}{ -u(qy +uw +tx)^{2}}{q^{4}t}U_{3}\genfrac{(}{)}{}{}{q^{4}tw}{(qy +uw +tx)^{3}}U_{23}\genfrac{(}{)}{}{}{ -tx}{qy +uw +tx}U_{24}\genfrac{(}{)}{}{}{ -z}{qy +uw +tx}$
\end{tabular}} \\
\hline
\end{tabular}\end{center}\par
\begin{center}
\begin{tabular}[c]{|c|}\hline
$x_{2} +x_{4} +x_{13}$ \\
\hline
$Z(X_{1} ,X_{3} ,\thinspace X_{4}X_{23} +X_{2}X_{34}) \cap V(X_{2} ,X_{4} ,\thinspace X_{2}X_{13} -X_{12}X_{23})$ \\
\hline
$m =(0 ,r ,0 ,t ,u ,v ,w ,x ,y ,z) ,\ tv +rw =0 ,\ r ,t , ,rx -uv \neq 0$ \\
\hline
{\begin{tabular}[c]{c}$T(\frac{t^{\frac{1}{3}}(rx -uv)^{\frac{2}{3}}}{r^{\frac{1}{3}}} ,1 ,\frac{1}{r} ,\frac{r^{\frac{2}{3}}t^{\frac{1}{3}}}{(rx -uv)^{\frac{1}{3}}})$$ \cdot $ \\
$U_{4}\genfrac{(}{)}{}{}{yu -rz}{t(rx -uv)}U_{3}\genfrac{(}{)}{}{}{ -vt^{\frac{1}{3}}r^{\frac{2}{3}}}{(rx -uv)^{\frac{1}{3}}}U_{1}\genfrac{(}{)}{}{}{u}{s^{\frac{2}{3}}t^{\frac{1}{3}}(rx -uv)^{\frac{2}{3}}}U_{23}\genfrac{(}{)}{}{}{(xy -vz)r^{\frac{5}{3}}}{t^{\frac{2}{3}}(rx -uv)^{\frac{4}{3}}}$
\end{tabular}} \\
\hline
\end{tabular}\end{center}\par
\begin{center}
\begin{tabular}[c]{|c|}\hline
$x_{2} +x_{4} +x_{34}$ \\
\hline
$Z(X_{1} ,X_{3} ,X_{2}X_{13} -X_{12}X_{23}) \cap V(X_{2} ,X_{4} ,X_{2}X_{34} +X_{4}X_{23})$ \\
\hline
$m =(0 ,r ,0 ,t ,u ,v ,w ,x ,y ,z) ,\ rx -uv =0 ,\ r ,t ,rw +tv \neq 0$ \\
\hline
{\begin{tabular}[c]{c}$T(\frac{rt^{3}}{(rw +tv)^{2}} ,\frac{rw +tv}{t} ,\frac{rw +tv}{rt} ,1)$ \\
$U_{3}\genfrac{(}{)}{}{}{ -tv}{rw +tv}U_{1}\genfrac{(}{)}{}{}{u(rw +tv)^{3}}{r^{2}t^{4}}U_{34}\genfrac{(}{)}{}{}{ -y}{rw +tv}U_{13}\genfrac{(}{)}{}{}{(rw +tv)^{2}(rz -yu)}{r^{2}t^{4}}$
\end{tabular}} \\
\hline
\end{tabular}\end{center}\par
\begin{center}
\begin{tabular}[c]{|c|}\hline
$x_{3} +x_{4}$ \\
\hline
$Z(X_{1} ,X_{2} ,X_{12}) \cap V(X_{3} ,X_{4})$ \\
\hline
$m =(0 ,0 ,s ,t ,0 ,v ,w ,x ,y ,z) ,\ s ,t \neq 0$ \\
\hline
$T(1 ,1 ,s^{\frac{2}{3}}t^{\frac{1}{3}} ,t^{\frac{1}{3}}s^{ -\frac{1}{3}})U_{4}\genfrac{(}{)}{}{}{ -w}{st}U_{2}\genfrac{(}{)}{}{}{t^{\frac{1}{3}}v}{s^{\frac{1}{3}}}U_{23}\genfrac{(}{)}{}{}{sy -vw}{s^{\frac{4}{3}}t^{\frac{2}{3}}}U_{12}\genfrac{(}{)}{}{}{t^{\frac{1}{3}}x}{s^{\frac{1}{3}}}U_{13}\genfrac{(}{)}{}{}{sz -wx}{s^{\frac{4}{3}}t^{\frac{2}{3}}}$ \\
\hline
\end{tabular}\end{center}\par
\begin{center}
\begin{tabular}[c]{|c|}\hline
$x_{1} +x_{2} +x_{34}$ \\
\hline
$Z(X_{3} ,X_{4}) \cap V(X_{1} ,X_{2} ,X_{34})$ \\
\hline
$m =(q ,r ,0 ,0 ,u ,v ,w ,x ,y ,z) ,\ q ,r ,w \neq 0$ \\
\hline
$T(q ,1 ,\frac{1}{r} ,\frac{wr^{2}}{q})U_{1}\genfrac{(}{)}{}{}{u}{qr}U_{3}\genfrac{(}{)}{}{}{ -vwr^{2}}{q}U_{23}\genfrac{(}{)}{}{}{rw(uv -rx)}{q^{2}}U_{34}\genfrac{(}{)}{}{}{ -y}{rw}U_{24}\genfrac{(}{)}{}{}{yu -rz}{qr^{2}w}$ \\
\hline
\end{tabular}\end{center}\par
\begin{center}
\begin{tabular}[c]{|c|}\hline
$x_{1} +x_{3} +x_{12} +x_{24}$ \\
\hline
$Z(X_{2} ,X_{4}) \cap V(X_{1} ,X_{3} ,\thinspace X_{1}X_{23} +X_{3}X_{12} ,\thinspace X_{3}X_{24} -X_{23}X_{34})$ \\
\hline
$m =(q ,0 ,s ,0 ,u ,v ,w ,x ,y ,z) ,\ \ q ,s ,qv +su ,sy -vw \neq 0$ \\
\hline
{\begin{tabular}[c]{c}$T(\frac{q^{\frac{2}{5}}(sy -vw)^{\frac{1}{5}}(qv +su)^{\frac{2}{5}}}{s^{\frac{2}{5}}} ,\frac{\thinspace (sy -vw)^{\frac{1}{5}}(qv +su)^{\frac{2}{5}}}{q^{\frac{3}{5}}s^{\frac{2}{5}}} ,\thinspace \frac{q^{\frac{2}{5}}s^{\frac{3}{5}}(sy -vw)^{\frac{1}{5}}}{(qv +su)^{\frac{3}{5}}} ,\frac{q^{\frac{2}{5}}(sy -vw)^{\frac{1}{5}}}{s^{\frac{2}{5}}(qv +su)^{\frac{3}{5}}})$$ \cdot $ \\
$U_{4}\genfrac{(}{)}{}{}{w(qv +su)}{q(vw -sy)}U_{2}\genfrac{(}{)}{}{}{qv}{qv +su}U_{23}\genfrac{(}{)}{}{}{ -x}{qv +su}U_{24}\genfrac{(}{)}{}{}{wx -sz}{q(sy -vw)}$
\end{tabular}} \\
\hline
\end{tabular}\end{center}\par
\begin{center}
\begin{tabular}[c]{|c|}\hline
$x_{2} +x_{3} +x_{14}$ \\
\hline
$Z(X_{1} ,X_{4}) \cap V(X_{2} ,X_{3} ,\thinspace X_{3}X_{12}X_{24} +X_{2}X_{34}X_{13} -X_{12}X_{23}X_{34} -X_{2}X_{3}X_{14})$ \\
\hline
$m =(0 ,r ,s ,0 ,u ,v ,w ,x ,y ,z) ,\ r ,s ,syu +rwx -uvw -rsz \neq 0$ \\
\hline
{\begin{tabular}[c]{c}$T(\frac{\sqrt{rsz +uvw -rwx -syu}}{r} ,r ,1 ,\frac{1}{s})U_{1}\genfrac{(}{)}{}{}{ru}{\sqrt{rsz +uvw -rwx -syu}}U_{2}\genfrac{(}{)}{}{}{v}{rs}$$ \cdot $ \\
$U_{4}\genfrac{(}{)}{}{}{ -sw}{\sqrt{rsz +uvw -rwx -syu}}U_{12}\genfrac{(}{)}{}{}{rx -uv}{s\sqrt{rsz +uvw -rwx -syu}}U_{34}\genfrac{(}{)}{}{}{vw -sy}{r\sqrt{rsz +uvw -rwx -syu}}$
\end{tabular}} \\
\hline
\end{tabular}\end{center}\par
\begin{center}
\begin{tabular}[c]{|c|}\hline
$x_{1} +x_{4} +x_{23}$ \\
\hline
$Z(X_{2} ,X_{3}) \cap V(X_{1} ,X_{4} ,X_{23})$ \\
\hline
$m =(q ,0 ,0 ,t ,u ,v ,w ,x ,y ,z) ,\ q ,t ,v \neq 0$ \\
\hline
{\begin{tabular}[c]{c}$T(qv ,v ,\frac{t}{qv^{2}} ,1)U_{2}\genfrac{(}{)}{}{}{ -tu}{q^{2}v^{3}}$$ \cdot $ \\
$U_{3}\genfrac{(}{)}{}{}{wqv^{2}}{t^{2}}U_{4}\genfrac{(}{)}{}{}{ -uw -tx -qy}{qtv}U_{23}\genfrac{(}{)}{}{}{ -x}{qv}U_{24}\genfrac{(}{)}{}{}{tx^{2} -qvz +qxy +uwx}{q^{2}tv^{2}}$
\end{tabular}} \\
\hline
\end{tabular}\end{center}\par
\begin{center}
\begin{tabular}[c]{|c|}\hline
$x_{2} +x_{4} +x_{34} +x_{13}$ \\
\hline
$Z(X_{1} ,X_{3}) \cap V(X_{2} ,X_{4} ,\thinspace X_{4}X_{23} +X_{2}X_{34} ,\thinspace X_{2}X_{13} -X_{12}X_{23})$ \\
\hline
$m =(0 ,r ,0 ,t ,u ,v ,w ,x ,y ,z) ,\ r ,t ,tv +rw ,rx -uv \neq 0$ \\
\hline
$$
{\begin{tabular}[c]{c}$T(\frac{t^{\frac{3}{5}}(rx -uv)^{\frac{4}{5}}}{r^{\frac{3}{5}}(rw +tv)^{\frac{2}{5}}} ,\frac{r^{\frac{2}{5}}(rw +tv)^{\frac{3}{5}}}{t^{\frac{2}{5}}(rx -uv)^{\frac{1}{5}}} ,\frac{(rw +tv)^{\frac{3}{5}}}{r^{\frac{3}{5}}t^{\frac{2}{5}}(rx -uv)^{\frac{1}{5}}} ,\frac{r^{\frac{2}{5}}t^{\frac{3}{5}}}{(rw +tv)^{\frac{2}{5}}(rx -uv)^{\frac{1}{5}}})$ \\
$U_{1}\genfrac{(}{)}{}{}{u(rw +tv)}{t(rx -uv)}U_{3}\genfrac{(}{)}{}{}{ -tv}{rw +tv}U_{23}\genfrac{(}{)}{}{}{y}{rw +tv}U_{13}\genfrac{(}{)}{}{}{rz -yu}{t(rx -uv)}$
\end{tabular}} \\
\hline
\end{tabular}\end{center}\par
\begin{center}
\begin{tabular}[c]{|c|}\hline
$x_{3} +x_{4} +x_{12}$ \\
\hline
$Z(X_{1} ,X_{2}) \cap V(X_{3} ,X_{4} ,X_{12})$ \\
\hline
$m =(0 ,0 ,s ,t ,u ,v ,w ,x ,y ,z) ,\ s ,t ,u \neq 0$ \\
\hline
$T(su ,\frac{t}{s^{2}u} ,s ,1)U_{4}\genfrac{(}{)}{}{}{ -w}{st}U_{2}\genfrac{(}{)}{}{}{s^{2}uv}{t}U_{23}\genfrac{(}{)}{}{}{su(sy -vw)}{t^{2}}U_{12}\genfrac{(}{)}{}{}{x}{su}U_{13}\genfrac{(}{)}{}{}{sz -wx}{s^{2}tu}$ \\
\hline
\end{tabular}\end{center}\par
\begin{center}
\begin{tabular}[c]{|c|}\hline
$x_{1} +x_{2} +x_{4}$ \\
\hline
$Z(X_{3} ,\thinspace X_{2}X_{34} +X_{4}X_{23}) \cap V(X_{1} ,X_{2} ,X_{4})$ \\
\hline
$m =(q ,r ,0 ,t ,u ,v ,w ,x ,y ,z) ,\ rw +tv =0 ,\ q ,r ,t \neq 0$ \\
\hline
$T(q ,1 ,\frac{1}{r} ,\sqrt{\frac{rt}{q}})U_{2}\genfrac{(}{)}{}{}{ -u}{qr}U_{3}\genfrac{(}{)}{}{}{ -v\sqrt{rt}}{\sqrt{q}}U_{23}\genfrac{(}{)}{}{}{ -x\sqrt{rt}}{\sqrt{q^{3}}}U_{34}\genfrac{(}{)}{}{}{uvt -rxt -qry}{^{\sqrt{q^{3}rt}}}U_{24}\genfrac{(}{)}{}{}{ -z\sqrt{r}}{\sqrt{q^{3}t}}$ \\
\hline
\end{tabular}\end{center}\par
\begin{center}
\begin{tabular}[c]{|c|}\hline
$x_{1} +x_{3} +x_{4}$ \\
\hline
$Z(X_{2} ,\thinspace X_{3}X_{12} +X_{1}X_{23}) \cap V(X_{1} ,X_{3} ,X_{4})$ \\
\hline
$m =(q ,0 ,s ,t ,u ,v ,w ,x ,y ,z) ,\ su +qv =0 ,\ q ,s ,t \neq 0$ \\
\hline
$T(\sqrt{\frac{qt}{s}} ,\sqrt{\frac{t}{qs}} ,s ,1)U_{3}\genfrac{(}{)}{}{}{w}{st}U_{2}\genfrac{(}{)}{}{}{v\sqrt{qs}}{\sqrt{t}}U_{23}\genfrac{(}{)}{}{}{y\sqrt{qs}}{\sqrt{t^{3}}}U_{12}\genfrac{(}{)}{}{}{sqy +stx -qvw}{\sqrt{qst^{3}}}U_{13}\genfrac{(}{)}{}{}{z\sqrt{s}}{\sqrt{qt^{3}}}$ \\
\hline
\end{tabular}\end{center}\par
\begin{center}
\begin{tabular}[c]{|c|}\hline
$x_{1} +x_{2} +x_{3}$\textsubscript {} \\
\hline
$Z(X_{4}) \cap V(X_{1} ,X_{2} ,X_{3})$ \\
\hline
$m =(q ,r ,s ,0 ,u ,v ,w ,x ,y ,z) ,\ q ,r ,s \neq 0$ \\
\hline
{\begin{tabular}[c]{c}$T(qrs ,rs ,s ,1)U_{1}\genfrac{(}{)}{}{}{u}{qr}U_{3}\genfrac{(}{)}{}{}{ -v}{rs}$$ \cdot $ \\
$U_{4}\genfrac{(}{)}{}{}{ -w}{qr^{2}s^{4}}U_{23}\genfrac{(}{)}{}{}{uv -rx}{qr^{2}s}U_{34}\genfrac{(}{)}{}{}{ -y}{qr^{3}s^{4}}U_{24}\genfrac{(}{)}{}{}{rwx -rsz +syu -uvw}{q^{2}r^{4}s^{5}}$
\end{tabular}} \\
\hline
\end{tabular}\end{center}\par
\begin{center}
\begin{tabular}[c]{|c|}\hline
$x_{1} +x_{2} +x_{4} +x_{34}$ \\
\hline
$Z(X_{3}) \cap V(X_{1} ,X_{2} ,X_{4} ,X_{2}X_{34} +X_{4}X_{23})$ \\
\hline
$m =(q ,r ,0 ,t ,u ,v ,w ,x ,y ,z) ,\ q ,r ,t ,rw +tv \neq 0$ \\
\hline
$$
{\begin{tabular}[c]{c}$T(\frac{q^{\frac{4}{5}}r^{\frac{1}{5}}(rw +tv)^{\frac{2}{5}}}{t^{\frac{1}{5}}} ,\frac{r^{\frac{1}{5}}(rw +tv)^{\frac{2}{5}}}{q^{\frac{1}{5}}t^{\frac{1}{5}}} ,\frac{(rw +tv)^{\frac{2}{5}}}{q^{\frac{1}{5}}r^{\frac{4}{5}}t^{\frac{1}{5}}} ,\frac{r^{\frac{1}{5}}t^{\frac{4}{5}}}{q^{\frac{1}{5}}(rw +tv)^{\frac{3}{5}}})$$ \cdot $ \\
$U_{2}\genfrac{(}{)}{}{}{ -u}{qr}U_{3}\genfrac{(}{)}{}{}{ -tv}{rw +tv}U_{23}\genfrac{(}{)}{}{}{ -tx}{q(rw +tv)}U_{34}\genfrac{(}{)}{}{}{ -qy -tx -uw}{q(rw +tv)}U_{24}\genfrac{(}{)}{}{}{ -z}{q(rw +tv}$
\end{tabular}} \\
\hline
\end{tabular}\end{center}\par
\begin{center}
\begin{tabular}[c]{|c|}\hline
$x_{1} +x_{3} +x_{4} +x_{12}$ \\
\hline
$Z(X_{2}) \cap V(X_{1} ,X_{3} ,X_{4} ,X_{1}X_{23} +X_{3}X_{12})$ \\
\hline
$m =(q ,0 ,s ,t ,u ,v ,w ,x ,y ,z) ,\ q ,s ,t ,qv +su \neq 0$ \\
\hline
{\begin{tabular}[c]{c}$T(\frac{q^{\frac{1}{5}}t^{\frac{1}{5}}(qv +su)^{\frac{3}{5}}}{s^{\frac{1}{5}}} ,\frac{t^{\frac{1}{5}}(qv +su)^{\frac{3}{5}}}{q^{\frac{4}{5}}s^{\frac{1}{5}}} ,\frac{q^{\frac{1}{5}}s^{\frac{4}{5}}t^{\frac{1}{5}}}{(qv +su)^{\frac{2}{5}}} ,\frac{q^{\frac{1}{5}}t^{\frac{1}{5}}}{s^{\frac{1}{5}}(qv +su)^{\frac{2}{5}}})$$ \cdot $ \\
$U_{3}\genfrac{(}{)}{}{}{w}{st}U_{2}\genfrac{(}{)}{}{}{qv}{qv +su}U_{23}\genfrac{(}{)}{}{}{qy}{t(qv +su)}U_{12}\genfrac{(}{)}{}{}{qy +tx +uw}{t(qv +su)}U_{13}\genfrac{(}{)}{}{}{z}{t(qv +su)}$
\end{tabular}} \\
\hline
\end{tabular}\end{center}\par
\begin{center}
\begin{tabular}[c]{|c|}\hline
$x_{2} +x_{3} +x_{4}$ \\
\hline
$Z(X_{1}) \cap V(X_{2} ,X_{3} ,X_{4})$ \\
\hline
$m =(0 ,r ,s ,t ,u ,v ,w ,x ,y ,z) ,\ r ,s ,t \neq 0$ \\
\hline
{\begin{tabular}[c]{c}$T(r^{3}s^{2}t ,1 ,\frac{1}{r} ,\frac{1}{rs})U_{4}\genfrac{(}{)}{}{}{ -w}{st}U_{2}\genfrac{(}{)}{}{}{v}{rs}$$ \cdot $ \\
$U_{1}\genfrac{(}{)}{}{}{u}{r^{4}s^{2}t}U_{23}\genfrac{(}{)}{}{}{sy -vw}{rs^{2}t}U_{12}\genfrac{(}{)}{}{}{x}{r^{4}s^{3}t}U_{13}\genfrac{(}{)}{}{}{rsz -rwx -syu +uvw}{r^{5}s^{4}t^{2}}$
\end{tabular}} \\
\hline
\end{tabular}\end{center}\par
\begin{center}
\begin{tabular}[c]{|c|}\hline
$x_{1} +x_{2} +x_{3} +x_{4}$ \\
\hline
$V(X_{1} ,X_{2} ,X_{3 ,}X_{4})$ \\
\hline
$m =(q ,r ,s ,t ,u ,v ,w ,x ,y ,z) ,\ q ,r ,s ,t \neq 0$ \\
\hline
{\begin{tabular}[c]{c}$T(q^{\frac{4}{5}}r^{\frac{3}{5}}s^{\frac{2}{5}}t^{\frac{1}{5}} ,\frac{r^{\frac{3}{5}}s^{\frac{2}{5}}t^{\frac{1}{5}}}{q^{\frac{1}{5}}} ,\frac{s^{\frac{2}{5}}t^{\frac{1}{5}}}{q^{\frac{1}{5}}r^{\frac{2}{5}}} ,\frac{t^{\frac{1}{5}}}{q^{\frac{1}{5}}r^{\frac{2}{5}}s^{\frac{3}{5}}})U_{2}\genfrac{(}{)}{}{}{ -u}{qr}U_{3}\genfrac{(}{)}{}{}{ -qv -su}{qrs}$$ \cdot $ \\
$U_{4}\genfrac{(}{)}{}{}{ -qrw -qtv -stu}{qrst}U_{12}\genfrac{(}{)}{}{}{x}{qrs}U_{34}\genfrac{(}{)}{}{}{ -qy -uw}{qrst}U_{24}\genfrac{(}{)}{}{}{qrwx -qrsz +qtvx +stux}{q^{2}r^{2}s^{2}t}$
\end{tabular}} \\
\hline
\end{tabular}\end{center}\par
We now verify that these $61$ orbits exhaust $\mathfrak{n} .$ Again, let $n =qx_{1} +rx_{2} +sx_{3} +tx_{4} +ux_{12} +vx_{23} +wx_{34} +xx_{13} +yx_{24} +zx_{14} =(q ,r ,s ,t ,u ,v ,w ,x ,y ,z)$ be an arbitrary element of $\mathfrak{n} .$ As in the above cases, we consider cases according to how many and which coordinates are $0.$  If $q ,r ,s ,t$ are all nonzero, then $n$ is in the regular orbit $B .(x_{1} +x_{2} +x_{3} +x_{4}) =V(X_{1} ,X_{2} ,X_{3} ,X_{4}) ,$ we we may assume at least one of these coordinates is zero.

First, consider the four cases when exactly one of these is $0.$  If $q =0 ,$ then $n \in Z(X_{1}) \cap V(X_{2} ,X_{3} ,X_{4}) =B .(x_{2} +x_{3} +x_{4}) .$  By symmetry (of the Dynkin diagram for $A_{4}) ,$ if $t =0$ then $n \in Z(X_{4}) \cap V(X_{1} ,X_{2} ,X_{3}) =B(x_{1} +x_{2} +x_{3}) .$  If $r =0 ,$ then $n \in Z(X_{2}) \cap V(X_{1} ,X_{3} ,X_{4}) =B .(x_{1} +x_{3} +x_{4}) \cup B(x_{1} +x_{3} +x_{4} +x_{12})$ - the union splits according to whether or not the value of $X_{1}X_{23} +X_{3}X_{12}$ is $0 ,$ so $n$ is in one of these two orbits.  By symmetry, if just $s =0 ,$ then $n \in Z(X_{3}) \cap V(X_{1} ,X_{2} ,X_{4}) =B .(x_{1} +x_{2} +x_{4}) \cup B .(x_{1} +x_{2} +x_{4} +x_{34}) .$

Next, we consider the six cases when exactly two of $q ,r ,s ,t$ are $0.$  The tables below handle these six cases, and the next four cases when exactly three of \thinspace $q ,r ,s ,t$ are $0.$  In the left column, we specify which coordinates are $0$ and the resulting locally closed set.  In the corresponding entry in the right column we list the orbit or union of orbits which comprise that locally closed set, showing in every case that $n$ belongs to one of the orbits:

\begin{center}
\begin{tabular}[c]{|c|c|}\hline
$0$ coordinates & Orbits \\
\hline
\hline
$q ,r$\par
$Z(X_{1} ,X_{2}) \cap V(X_{3} ,X_{4})$ & 
{\begin{tabular}[c]{c}$B .(x_{3} +x_{4} +x_{12}) \cup B .(x_{3} +x_{4})$ \\
according to the value of $X_{12}$
\end{tabular}}\par
\\
\hline
$q ,s$\par
$Z(X_{1} ,X_{3}) \cap V(X_{2} ,X_{4})$
& 
{\begin{tabular}[c]{c}$B .(x_{2} +x_{4} +x_{34} +x_{123}) \cup B .(x_{2} +x_{4} +x_{13})$ \\
$ \cup B .(x_{2} +x_{4} +x_{34}) \cup B .(x_{2} +x_{4})$ \\
according to the values of \\
$X_{4}X_{23} +X_{2}X_{34}$ and $X_{2}X_{13} -X_{12}X_{23}$
\end{tabular}} \\
\hline
$q ,t$\par
$Z(X_{1} ,X_{4}) \cap V(X_{2} ,X_{3})$
& 
{\begin{tabular}[c]{c}$B .(x_{2} +x_{3} +x_{14}) \cup B .(x_{2} +x_{3})$ \\
according to the value of  \\
$X_{3}X_{12}X_{24} +X_{2}X_{34}X_{13} -X_{12}X_{23}X_{34} -X_{2}X_{3}X_{14}$
\end{tabular}} \\
\hline
$r ,s$\par
$Z(X_{2} ,X_{3}) \cap V(X_{1} ,X_{4})$
& 
{\begin{tabular}[c]{c}$B .(x_{1} +x_{4} +x_{23}) \cup B .(x_{1} +x_{4} +x_{24}) \cup B .(x_{1} +x_{4})$ \\
according to the values of \\
$X_{23}$ and $X_{1}X_{24} +X_{12}X_{34} +X_{4}X_{13}$
\end{tabular}} \\
\hline
$r ,t$\par
$Z(X_{2} ,X_{4}) \cap V(X_{1} ,X_{3})$
& 
{\begin{tabular}[c]{c}$B .(x_{1} +x_{3} +x_{12} +x_{24}) \cup B .(x_{1} +x_{3} +x_{12})$ \\
$ \cup B .(x_{1} +x_{3} +x_{24}) \cup B .(x_{1} +x_{3})$ \\
according to the values of  \\
$X_{1}X_{23} +X_{3}X_{12}$ and $X_{3}X_{24} -X_{23}X_{34}$
\end{tabular}} \\
\hline
$s ,t$\par
$Z(X_{3} ,X_{4}) \cap V(X_{1} ,X_{2})$
& 
{\begin{tabular}[c]{c}$B .(x_{1} +x_{2} +x_{34}) \cup B .(x_{1} +x_{2})$ \\
according to the values of $X_{34}$
\end{tabular}} \\
\hline
\multicolumn{2}{|c|}{Orbits when exactly $2$ of $q ,r ,s ,t$ are $0$} \\
\hline
\end{tabular}\end{center}\par
\begin{center}
\begin{tabular}[c]{|c|c|}\hline
$30 +x$ coordinates & Orbits \\
\hline
\hline
$q ,r ,s$\par
$Z(X_{1} ,X_{2} ,X_{3}) \cap V(X_{4})$
& 
{\begin{tabular}[c]{c}$B .(x_{4} +x_{12} +x_{23}) \cup (x_{4} +x_{12} +x_{13})$ \\
$ \cup B .(X_{4} +x_{23}) \cup B .(x_{4} +x_{12}) \cup B .x_{4}$ \\
according to the values of  \\
$X_{12} ,X_{23} ,X_{13}$ and $X_{4}X_{13} +X_{12}X_{34}$
\end{tabular}}\par
\\
\hline
$q ,r ,t$\par
$Z(X_{1} ,X_{2} ,X_{4}) \cap V(X_{3})$
& 
{\begin{tabular}[c]{c}$B .(x_{3} +x_{12} +x_{24}) \cup B .(x_{3} +x_{12})$ \\
$ \cup B .(x_{3} +x_{24}) \cup B .(x_{3} +x_{14}) \cup B .x_{3}$ \\
according to the values of \\
$X_{12} ,\thinspace X_{3}X_{24} -X_{34}X_{34}$ and $X_{3}X_{14} -X_{34}X_{13}$
\end{tabular}}\par
\\
\hline
$q ,s ,t$\par
$Z(X_{1} ,X_{3} ,X_{4}) \cap V(X_{2})$
& 
{\begin{tabular}[c]{c}$B .(x_{2} +x_{34} +x_{13}) \cup B .(x_{2} +x_{13})$ \\
$ \cup B .(x_{2} +x_{34}) \cup B .(x_{2} +x_{14}) \cup B .x_{2}$ \\
according to the values of \\
$X_{34} ,\thinspace X_{2}X_{13} -X_{12}X_{23}$ and $X_{2}X_{14} -X_{12}X_{24}$
\end{tabular}}\par
\\
\hline
$r ,s ,t$\par
$Z(X_{2} ,X_{3} ,X_{4}) \cap V(X_{1})$
& 
{\begin{tabular}[c]{c}$B .(x_{1} +x_{23} +x_{34}) \cup B .(x_{1} +x_{34} +x_{24})$ \\
$ \cup B .(x_{1} +x_{23}) \cup B .(x_{1} +x_{34}) \cup B .x_{1}$ \\
according to the values of \\
$X_{23} ,X_{34} ,X_{24}$ and $X_{1}X_{24} +_{12}X_{34}$
\end{tabular}} \\
\hline
\multicolumn{2}{|c|}{Orbits when exactly $3$ of $q ,r ,s ,t$ are $0$} \\
\hline
\end{tabular}\end{center}\par
The last case, when all four of $q ,r ,s ,t$ are zero, breaks into several subcases, according to how many of $u ,v ,w$ are also \thinspace $0.$  First, if none of these three are $0 ,$ then $n \in Z(X_{1} ,X_{2} ,X_{3} ,X_{4}) \cap V(X_{12} ,X_{23} ,X_{34}) =B .(x_{12} +x_{23} +x_{34}) .$  Next, consider the three cases when exactly one of these is $0 :$

\begin{center}
\begin{tabular}[c]{|c|c|}\hline
(height $2)$ coordinates with value $0$ & Orbits \\
\hline
\hline
$u$\par
$Z(X_{1} ,X_{2} ,X_{3} ,X_{4} ,X_{12}) \cap V(X_{23} ,X_{34})$
& $B .(x_{23} +x_{34})$ \\
\hline
$v$\par $Z(X_{1} ,X_{2} ,X_{3} ,X_{4} ,X_{23}) \cap V(X_{12} ,X_{34})$
& $B .(x_{12} +x_{34})$ \\
\hline
$w$\par $Z(X_{1} ,X_{2} ,X_{3} ,X_{4} ,X_{34}) \cap V(X_{12} ,X_{23})$
& $B .(x_{12} +x_{23})$ \\
\hline
\multicolumn{2}{|c|}{Orbits when exactly $1$ of $u ,v ,w$ are $0$} \\
\hline
\end{tabular}\end{center}\par
Now, the three cases when exactly two of these are zero:

\begin{center}
\begin{tabular}[c]{|c|c|}\hline
(height $2)$ coordinates with value $0$ & Orbits \\
\hline
\hline
$u ,v$\par
$Z(X_{1} ,X_{2} ,X_{3} ,X_{4} ,X_{12} ,X_{23}) \cap V(X_{34})$
& 
{\begin{tabular}[c]{c}$B .(x_{34} +x_{13}) \cup B .x_{34}$ \\
according to the value of  \\
$X_{13}$ 
\end{tabular}} \\
\hline
$u ,w$\par $Z(X_{1} ,X_{2} ,X_{3} ,X_{4} ,X_{12} ,X_{34}) \cap V(X_{23})$
& 
{\begin{tabular}[c]{c}$B .(x_{23} +x_{14}) \cup B .x_{23}$ \\
according to the value of  \\
$X_{13}X_{24} -X_{23}X_{14}$ 
\end{tabular}} \\
\hline
$v ,w$\par $Z(X_{1} ,X_{2} ,X_{3} ,X_{4} ,X_{23} ,X_{34}) \cap V(X_{12})$
& 
{\begin{tabular}[c]{c}$B .(x_{12} +x_{24}) \cup B .x_{12}$ \\
according to the value of  \\
$X_{24}$
\end{tabular}} \\
\hline
\multicolumn{2}{|c|}{Orbits when exactly $2$ of $u ,v ,w$ are $0$} \\
\hline
\end{tabular}\end{center}\par
That leaves the case when $q ,r ,s ,t ,u ,v ,w$ are all $0 ,$ which consists of further subcases.  If \thinspace $x ,y$ are both nonzero (coordinates of root vectors for height $3$ roots), then $n \in Z(X_{1} ,X_{2} ,X_{3} ,X_{4} ,X_{12} ,X_{23} ,X_{34}) \cap V(X_{13} ,X_{24}) =B .(x_{13} +x_{24}) .$ The remaining cases are tabulated below:

\begin{center}
\begin{tabular}[c]{|c|c|}\hline
(height $3)$ coordinates with value $0$ & Orbits \\
\hline
\hline
$x$\par
$Z(X_{1} ,X_{2} ,X_{3} ,X_{4} ,X_{12} ,X_{23} ,X_{34} ,X_{13}) \cap V(X_{24})$
& $B .x_{24}$ \\
\hline
$y$\par $Z(X_{1} ,X_{2} ,X_{3} ,X_{4} ,X_{12} ,X_{23} ,X_{34} ,X_{24}) \cap V(X_{13})$
& $B .x_{13}$
\\
\hline
$x ,y$\par $Z(X_{1} ,X_{2} ,X_{3} ,X_{4} ,X_{12} ,X_{23} ,X_{34} ,X_{13} ,X_{24})$
&
{\begin{tabular}[c]{c}$B .x_{14} \cup B.0$ \\
according to the value of \\
$X_{14}$
\end{tabular}} \\
\hline
\multicolumn{2}{|c|}{all remaining cases.} \\
\hline
\end{tabular}\end{center}\par
Thus, we have shown that the $61$ given orbits exhaust all of $\mathfrak{n} .$ It remains to check the containments indicated in the Hasse diagram of the closure order.  But they follow exactly as in the previous case. \bigskip 
\end{proof}

Since the Hasse diagram for the case $A_{4}$ is so cluttered and hard to read, we provide an enlarged version of it below.:

\includegraphics[ width=5.531250000000001in,]{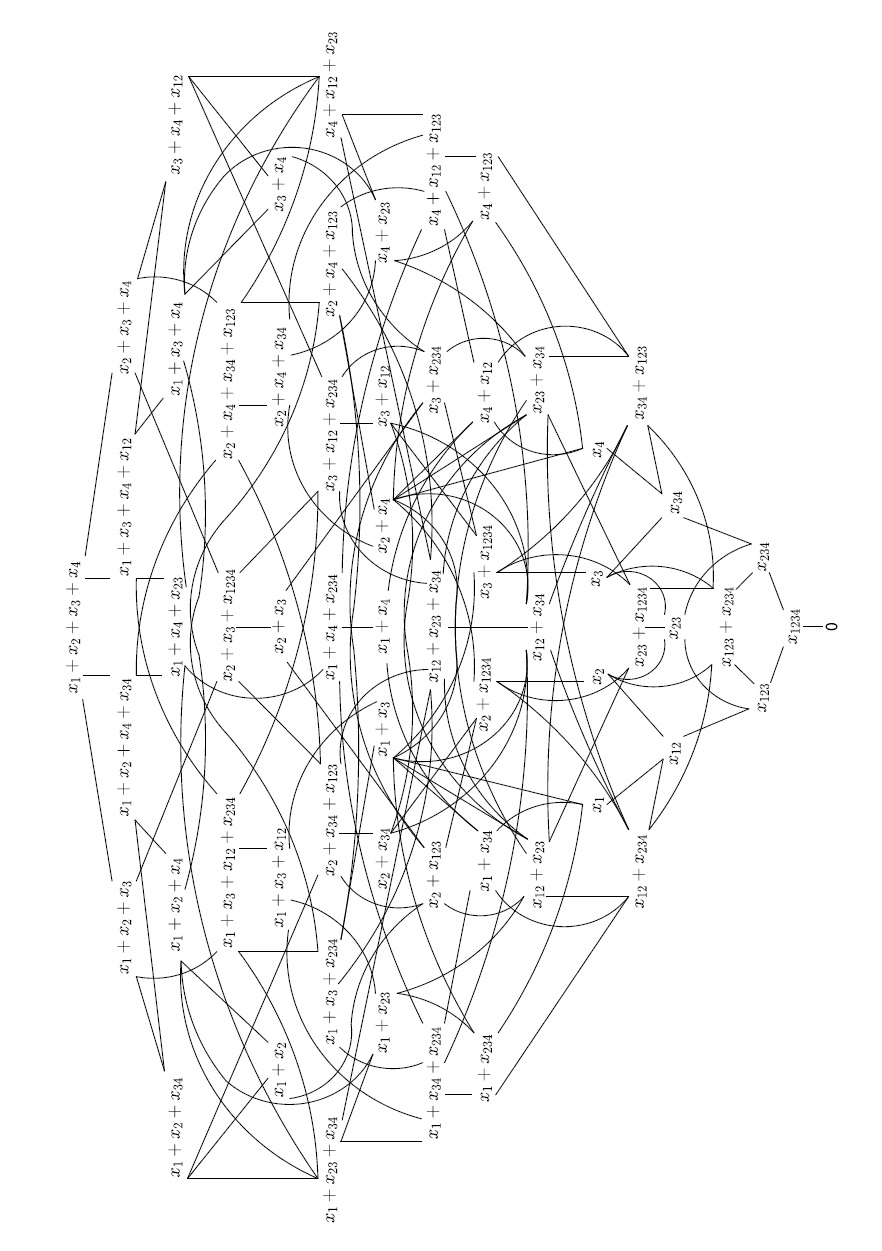}

\begin{center}{\Large References }\bigskip \end{center}\par

[BHRZ] T. Brustle, L. Hille, G. Rohrle, and G, Zwara, \textit{The Bruhat-Chevalley
Order of Parabolic Group Actions in General Linear Groups and Degeneration for }
$\Delta $
\textit{-Filtered Modules}, Advances in Mathematics, \textbf{148}, pp.
203-242, 1999. \bigskip

[Bou] N. Bourbaki, \textit{Elements of Mathematics, Lie Groups and Lie Algebras,
Chapters}\textit{4-6}, Springer, 2002 \bigskip

[BH] H. B{\"u}rgstein and W. H. Hesselink, \textit{Algorithmic Orbit Classication for
Some Borel Group Actions}, Composito Mathematica\textbf{ 61}, pp 3-41, 1987 \bigskip

[BV] M. Burkhart and D. Vella, \textit{Nilpotent Orbits for
Borel Subgroups of $SO_{5}(k)$, }unpublished. \bigskip

 [C] R. W. Carter, \textit{Finite Groups of Lie Type}, Wiley Classics
Library, 1993\bigskip

[CM] D. H. Collingwood and W. M. McGovern, \textit{Nilpotent Orbits in Semisimple
Lie Algebras}, Van Nostrand Rheinhold Mathematics Series, 1993 \bigskip

[H1] J. E. Humphreys,\textit{ Introduction to Lie Algebras and Representation
Theory}, Springer, Graduate Texts in Mathematics \textbf{9}, 1972 \bigskip

[H2] J.E. Humphreys,\textit{ Linear Algebraic Groups}, Springer, Graduate
Texts in Mathematics \textbf{21}, 1975 \bigskip

 [HR] L. Hille and G. R{\"o}hrle\textit{, A Classification of Parabolic Subgroups of
Classical Groups with a Finite Number of Orbits on the Unipotent Radical}, Transformation
Groups, Vol. 4, No. 1, pp. 35-52, 1999 \bigskip

 [J] J. Jantzen, \textit{Nilpotent Orbits in Representation Theory}, in
Lie Theory - Lie Algebras and Representations, Springer Progress in Mathematics \textbf{228,}
pp. 1-211, 2004. \bigskip

 [K] V. V. Kashin, \textit{Orbits of adjoint and coadjoint actions of Borel
subgroups of semisimple algebraic groups}. (In Russian) Questions of Group Theory and
Homological Algebra, Yaroslavl, 141-159, 1990 \bigskip

\bigskip [P] V. L. Popov, \textit{A Finiteness Theorem
for Parabolic Subgroups of Fixed Modality}, Indag. Mathem., N.S. \textbf{8} (1),
pp.125-132, 1997 \bigskip

[PR] V. L. Popov and G. R{\"o}hrle, \textit{On the Number of Orbits of a Parabolic
Subgroup on its Unipotent Radical}, Algebraic Groups and Lie Groups, G. I. Lehrer, Ed.,
Australian Mathematical Society, Vol. \textbf{9}, Cambridge University Press, 1997 \bigskip

 [R1] G. R{\"o}hrle, \textit{Parabolic Subgroups of Positive Modality},
Geometriae Dedicata,\textbf{ 60}, pp.163-186, 1996 \bigskip

[R2] G. R{\"o}hrle,
\textit{On the Modality of Parabolic
Subgroups of Linear Algebraic Groups}, Manuscripta Math. 98, pp, 9-20, 1999 \bigskip

\begin{center}
\end{center}\par

\end{document}